\documentclass[12pt,a4paper]{article}
\usepackage[utf8]{inputenc}
\usepackage[T1]{fontenc}
\usepackage{amsmath}
\usepackage{amsfonts}
\usepackage{amsthm}
\usepackage{color}

\usepackage[margin=1in]{geometry}

\newcommand{\cL}{{\mathbb {L}}}

\newcommand{\bpf}{\begin{preuve}}
\newcommand{\epf}{ \end{preuve} \medskip}

\newcommand{\benum}{\begin{enumerate}}
\newcommand{\eenum}{\end{enumerate}}

\newcommand{\bitem}{\begin{itemize}}
\newcommand{\eitem}{\end{itemize}}

\newcommand{\brmq}{\begin{rmq}}
\newcommand{\ermq}{\end{rmq}}

\newcommand{\brmqs}{\begin{rmqs}}
\newcommand{\ermqs}{\end{rmqs}}

\newcommand{\bapp}{\begin{application}}
\newcommand{\eapp}{\end{application}}

\newcommand{\bapps}{\begin{applications}}
\newcommand{\eapps}{\end{applications}}

\newcommand{\bdefi}{\begin{definition}}
\newcommand{\edefi}{\end{definition}}

\newcommand{\beq}{\begin{equation}}
\newcommand{\eeq}{\end{equation}}

\def\bpm{\begin{pmatrix}}
\def\epm{\end{pmatrix}}

\newcommand{\bcas}{\begin{cases}}
\newcommand{\ecas}{\end{cases}}

\newcommand{\bex}{\begin{exemp}}
\newcommand{\eex}{\end{exemp}}

\newcommand{\bexs}{\begin{exemps}}
\newcommand{\eexs}{\end{exemps}}

\newcommand{\bprop}{\begin{proposition}}
\newcommand{\eprop}{\end{proposition}}

\newcommand{\bthm}{\begin{theoreme}}
\newcommand{\ethm}{\end{theoreme}}

\newcommand{\bcor}{\begin{corollaire}}
\newcommand{\ecor}{\end{corollaire}}

\newcommand{\blem}{\begin{lemme}}
\newcommand{\elem}{\end{lemme}}

\newcommand{\beqna}{\begin{eqnarray}}
\newcommand{\eeqna}{\end{eqnarray}}

\newcommand{\beqnas}{\begin{eqnarray*}}
\newcommand{\eeqnas}{\end{eqnarray*}}

\newcommand{\cA}{{\mathcal A}}

\newcommand{\Cov}{{\rm Cov}} 

\def\tr{\textmd{trace}\,}

\def\Lip{\textsc{Lipschitz}}

\def\Id{{\rm{Id}}} 

\def\cA{{\mathcal A }}
\def\cB{{\mathcal B }}
\def\cC{{\mathcal C}}

\def\cH{{\mathcal  H}}
\def\cI{{\mathcal  I}}
\def\cK{{\mathcal  K}}
\def\cL{{\mathcal L }}

\def\cN{{\mathcal N }}
\def\cP{{\mathcal P }}

\def\cV{{\mathcal V}}

\def\bbE{{\mathbb{E}}}

\newcommand{\bbN}{{\mathbb {N}}}
\newcommand{\bbP}{{\mathbb P}}
\newcommand{\bbR}{{\mathbb {R}}}
\newcommand{\bbS}{{\mathbb {S}}} 

\newcommand{\bbZ}{\mathbb {Z}}

\def\un{{\mathbf{1}}}
\def\bfA{{\mathbf{A}}}

\def\bfe{{\mathbf{e}}}

\def\bfq{{\mathbf{q}}}
\def\bfQ{{\mathbf{Q}}}

\def\bfu{{\mathbf{u}}}
\def\bfV{{\mathbf{V}}}
\def\bfv{{\mathbf{v}}}
\def\bfW{{\mathbf{W}}}
\def\bfw{{\mathbf{w}}}
\def\bfX{{\mathbf{X}}}
\def\bfx{{\mathbf{x}}}
\def\bfY{{\mathbf{Y}}}

\def\bfZ{{\mathbf{Z}}}
\def\bfz{{\mathbf{z}}}
\def\bf\Sigma{{\mathbf{\Sigma}}}
\def\bdeta{{\boldsymbol{\eta}}}
\def\bdalpha{{\boldsymbol{\alpha}}}
\def\bfomega{{\boldsymbol{\omega}}}

\newtheorem{theoreme}{Theorem}[section]
\newtheorem{lemme}[theoreme]{Lemma}
\newtheorem{definition}[theoreme]{Definition}
\newtheorem{proposition}[theoreme]{Proposition}
\newtheorem{corollaire}[theoreme]{Corollary}

\newenvironment{exemp}{\noindent{\bf Example. --- }}{\par}
\newenvironment{exemps}{\noindent{\bf Examples}\benum}{\eenum\par}
\newtheorem{rmq}[theoreme]{Remark}
\newtheorem{rmqs}[theoreme]{Remarks}

\newenvironment{preuve}{\noindent{\it Proof. --- }}
{\hfill\rule{1.3mm}{2mm}\par}
\newenvironment{application}{\noindent{\bf Application. --- }}{\par}
\newenvironment{applications}{\noindent{\bf Applications. ---
}\benum}{\eenum\par}

\theoremstyle{definition}

\title{The convergence rates for the superdiffusion 
in the Boltzmann-Grad limit of the periodic Lorentz gas}

\author{Songzi Li\thanks{School of Mathematics, Renmin University of China, 59, Zhongguancun Da Jie, Beijing, 100872, China ({\sf sli@ruc.edu.cn}). Supported by NSFC No.~11901569 and by fund (No.2018030249) from Renmin University of China.}}
\date{\today}

\begin{document}
\maketitle

\begin{abstract}
In this article, we obtain the rates of convergence for superdiffusion in the Boltzmann-Grad limit of the periodic Lorentz gas, which is one of the fundamental models to study diffusions in deterministic systems. In their seminal work, Marklof and Str\"ombergsson proved the Boltzmann-Grad limit of the periodic Lorentz gas~\cite{M-Sannals2}, and then Marklof and T\'oth established a superdiffusive central limit theorem in large time for the Boltzmann-Grad limit~\cite{M-T16}. Based on their work, we apply Stein's method to derive the convergence rates for the superdiffusion in the Boltzmann-Grad limit of the periodic Lorentz gas. For the discrete time displacement the rate of convergence in Wasserstein distance is obtained, while in the context of the continuous time displacement the result is presented for the Berry-Essen type bound. 
\end{abstract}



\section{Introduction and main results}


Since its introduction~\cite{Lor1905} in 1905 to study the motion of electrons in metals, the Lorentz gas has been considered as a basic model to study diffusions in deterministic systems and to understand Brownian motion from microscopic perspectives. It is a dynamical system corresponding with a tagged particle moving in an infinitely array of fixed scatterers. There have been many research on this topic, involving dynamical systems, probability, differential equations, etc. 

The Boltzmann-Grad limit of the Lorentz gas corresponds to the dynamics of the tagged particle when taking the radius of the scatterers to zero. It was conjectured by Lorentz~\cite{Lor1905} that the limiting system is governed by a linear Boltzmann equation. Indeed it is true when the configuration of the scatterers is random, see Gallavotti~\cite{Gal1969},  Spohn~\cite{Spohn78}, C. Boldrighini, L.A. Bunimovich and Y.G. Sinai~\cite{BBS1983}. As for the periodic Lorentz gas, Golse~\cite{Gol06} pointed out that the linear Boltzmann equation fails. Later Boca-Zaharescu~\cite{Boca-Zaha07} gave the explicit formula of the limiting distribution of the free path length for the two-dimensional periodic Lorentz gas, see also Caglioti-Golse~\cite{Cag-Golse03},  \cite{Cag-Golse08}.  As for the high dimension, Marklof-Str\"ombergesson proved that the Boltzmann-Grad limit  of the periodic Lorentz gas exists and gave its descriptions in their seminal papers \cite{M-Sannals1}, \cite{M-Sannals2} and \cite{M-Sgafa}. 

When the radius of the scatterers is fixed, the diffusive property of the Lorentz gas has been extensively studied.  The classical result was derived by Bunimovich-Sinai \cite{Buni-Sinai1981}, where they proved the central limit theorem for the planar periodic Lorentz gas with finite horizon. Bleher~\cite{Bleh92} was the first to point out the super-diffusive nature of the periodic Lorentz gas with infinite horizon. This fact was later confirmed by Sz\'asz-Varju in \cite{Sz-Var07} for the planar billiard map, and by Chernov-Dolgopyat in \cite{Dol-Chern09} for the planar billiard flow.  

In~\cite{M-T16}, Marklof and T\'oth studied the long-time behavior of the Boltzmann-Grad limit of the periodic Lorentz gas by proving a superdiffusive central limit theorem. To be precise, let $\bfX_{t}$ be the continuous time displacement of the Boltzmann-Grad limit of the periodic Lorentz gas, and $\bfQ_{n}$ be its discrete time displacement. Then under assumptions $\bfA$ (which we will explain in Section 2), the following superdiffusive central limit theorem holds.  

\bthm\label{superdiffusion}
Let $d \geq 2$. Under assumptions $\bfA$ on the initial data, we have 
\beqnas
\frac{\bfX_{t}}{\Sigma_{d}\sqrt{t \log t}} \Rightarrow \cN(0, \Id_{d}),
\eeqnas
and
\beqnas
\frac{\bfQ_{n}}{\sigma_{d}\sqrt{n \log n}} \Rightarrow \cN(0, \Id_{d}).
\eeqnas
 where $\sigma^{2}_{d} =  \frac{2^{2-d}}{2d^{2}(d+1)\zeta(d)}$,  $\Sigma^{2}_{d} =  \frac{2^{2-d}}{2d^{2}(d+1)\bar{\xi}\zeta(d)}$ and $\bar{\xi}$ is the mean free flight length. 
\ethm

As the above result shows,  the Boltzmann-Grad limit of the periodic Lorentz gas is a superdiffusion, in that the scaling factor is now $\frac{1}{\sqrt{n \log n}}$ instead of $\frac{1}{\sqrt{n}}$ in the classical central limit theorem. Then together with the Boltzmann-Grad limit of the periodic Lorentz gas~\cite{M-Sannals2}, it was derived in~\cite{M-T16} that the Brownian motion can been realized as the limit of the periodic Lorentz gas, first under the Boltzmann-Grad limit and then the large time limit. 

In this paper, our aim is to describe how fast the distribution of the Boltzmann-Grad limit of the periodic Lorentz gas,  converges to the Gaussian distribution. More precisely, we first show the explicit convergence rates of $\bfQ_{n}$ in terms of the Wasserstein distance. Recall that the Wasserstein distance $d_{W_{1}}$ between two random variables $\bfX$ and $\bfY$ is defined by
$$
d_{W_{1}}(\bfX, \bfY) = \sup_{h \in \Lip(1)}|\bbE h(\bfX) - \bbE h(\bfY)|.
$$

Under the same assumptions  $\bfA$ as in Theorem~\ref{superdiffusion} (which will be explained in Section 2), we have the following estimates hold.
\bthm\label{thm1}
Let $\bfZ$ be standard Gaussian random variable on $\bbR^{d}$.  Under the assumptions $\bfA$, we have when $d=2$,
\beqna\label{thm1.est1}
d_{W_{1}}(\frac{\bfQ_{n}}{\sigma_{d}\sqrt{n \log n}}, \bfZ) \leq  O(\sqrt{\frac{\log \log n}{\log n}}),
 \eeqna
 and  when $d \geq 3$
 \beqna\label{thm1.est2}
d_{W_{1}}(\frac{\bfQ_{n}}{\sigma_{d}\sqrt{n \log n}}, \bfZ) \leq  O(\frac{1}{\sqrt{\log n}}).
 \eeqna

\ethm

As for the continuous time displacement $\bfX_{t}$, the convergence rates for the Berry-Essen bound are obtained. 

\bthm\label{thm3}
Let $\bfz \in \bbR^{d}$. Under the assumptions $\bfA$,  we have when $d =2$,  
\beqna\label{thm3.cont1}
|\bbP(\frac{\bfX_{t}}{\Sigma_{d}\sqrt{t \log t}} \leq \bfz)  - \Phi(\bfz) | \leq O(\frac{(\log \log t)^{\frac{1}{4}}}{( \log t)^{\frac{1}{4}}}),
\eeqna
and when $d \geq 3$,
\beqna\label{thm3.cont2}
|\bbP(\frac{\bfX_{t}}{\Sigma_{d}\sqrt{t \log t}} \leq \bfz)  - \Phi(\bfz) | \leq O(\frac{1}{( \log t)^{\frac{1}{4}}}),
\eeqna
where $\Phi(\bfz)$ is the Gaussian distribution function on $\bbR^{d}$.
\ethm
 
The above result shows that the convergence rate for the superdiffusion in the Boltzmann-Grad limit of the periodic Lorentz gas are quite slow compared to that of the classical central limit theorem, which is $O(\frac{1}{\sqrt{n}})$. 

The general method to obtain the convergence rate in the central limit theorem is characteristic function (Fourier analysis). 
However,  the displacement  $\bfQ_{n}$ (or $\bfX_{t}$) as a process is very complicated: it involves many random variables which depend on each other, and the transition kernel is hard to deal with, which makes the method of  characteristic function impossible to work. 

The main tool of our proof is the Stein's method. Starting from a basic observation of the Gaussian distribution, Stein's method for normal approximation was first introduced by Stein \cite{Stein72}.  Since then, it has been shown to be powerful in many mathematical fields, for that the ideas of Stein's method is abstract enough to be extended to other distribution approximation (such as Poisson distribution) and even diffusions, and that it also works for a variety of metrics. There are lots of  comprehensive study on Stein's method, see for example \cite{C-S11},  \cite{notes_cha} and \cite{ross11}. Generally Stein's method is separated into two steps. The first one is through the so-called Stein's equation, which turns the difference between two distributions into the expectation with respect to the aimed distribution of a certain space of functions. The second step is to estimate the expectation derived in the first step. 

 The advantage of Stein's method is that it allows us to localize the difference between the distributions, which is efficient for complicated situations as $\bfQ_{n}$ (the displacement of the Boltzmann-Grad limit of the periodic Lorentz gas), such that with the moment estimates we are able to bound each term and derive the final results.  To apply Stein's method to our case, the first step is quite straightforward, since we are in the case of normal approximation. For the second step, inspired by the work of Chatterjee-Meckes~\cite{C-M07},  we manage to adapt the idea of exchangeable pairs for multivariate normal approximation to our case.  The technique of "exchangeable pairs" was first introduced by Stein \cite{Stein72} and given in details for the univariate case in \cite{Stein86}. The multivariate version of exchangeable pairs was first studied by Chatterjee-Meckes~\cite{C-M07} and then Reinert-R\"ollin \cite{RR09}. By applying the method of "exchangeable pairs" to our case, we directly obtain an explicit error.

To obtain the final estimates, it still requires bounding some specific forms of expectations in the error term, which are analogous to the moment bounds when applying Stein's method in the classical settings. Recall that the key point in the proof of superdiffusive central limit theorem in \cite{M-T16} is the spectral gap of the transition kernel, which also plays an important role in our work. Thanks to it, we are able to control the moments to get an explicit convergence rate. We point out that  for the continuous time displacement $\bfX_{t}$,  stronger estimates are required than those in the discrete case. This is the reason why the convergence rate is only shown for the Berry-Essen type bounds. 

As far as we know, there are few references regarding the convergence rates of Lorentz gas, or of superdiffusions. We mention the work of F. P\'ene~\cite{Pene}, where she proved the convergence rate $O(\frac{1}{\sqrt{n}})$ for the central limit theorem of the Sinai billiard,  in the sense of Kantorovich distance. 


As the referee pointed out, it would be also interesting if one considers the low density limit and the large time limit together, i.e., taking a joint limit $r \rightarrow 0$ and $t \rightarrow \infty$. Recall that in dimension two, Chernov-Dolgopyat~\cite{Dol-Chern09} derived a superdiffusive central limit theorem and invariance principle for the billiard flow when the radius is fixed. More recently, Lutsko-T\'oth~\cite{lutsko-toth2020} proved the invariance principle for a random Lorentz gas under the Boltzmann-Grad limit and large time limit simultaneously. Therefore to extend such results to the periodic Lorentz gas in any dimension is quite natural and meaningful . 


%

This paper is organized as follows. In Section 2 we offer some preliminaries about this topic, including the descriptions of the Boltzmann-Grad limit of periodic Lorentz gas, as well as some fundamental facts on Stein's method and exchangeable pairs. In Section 3 we prove the key estimates which are needed to bound the error term appearing in the Stein's method. In Section 4 and 5 we give the proofs of Theorem~\ref{thm1} and Theorem~\ref{thm3} respectively .

\section{Preliminaries}

\subsection{ The Boltzmann-Grad limit of the periodic Lorentz gas}
In the first part of this section we provide a more detailed description of the Boltzmann-Grad limit of the periodic Lorentz gas. Here we use the same notations as in \cite{M-T16}. For a comprehensive study, see Marklof-Str\"ombergsson~\cite{M-Sannals1, M-Sannals2}, Marklof-Toth~\cite{M-T16},  Marklof~\cite{M2014} and references therein.

We start with the settings of the periodic Lorentz gas. Let $\cL \subset \bbR^{d}$ be a fixed Euclidean lattice of covolume one. At each point of $\cL$ there is a sphere of radius $r$ as an obstacle. Define $\cL_{r} = r^{\frac{d-1}{d}}\cL$. We consider a test particle that moves at the speed $\bfv$ with $\|\bfv\| =1$ in the space $\cK_{r} = \bbR^{d} \backslash (\cL_{r} + r\cB^{d}_{1})$. Assume that there is no external potential such that the particle always moves in straight lines until it hits an obstacle, and we also assume that all the collisions are elastic. 

Denote by $\bfq_{n} = \bfq_{n}(\bfq_{0}, \bfv_{0}) \in \partial \cK_{r}$, the location where the test particle with initial condition $(\bfq_{0}, \bfv_{0})$ leaves the nth scatterer. Assume that $\bfq_{0}  \in \partial \cK_{r}$, such that the distribution of $\bfq_{0}$ can be obtained from $\bfv_{0}$. For the continuous time version, define $\bfx_{t} = \bfx_{t}(\bfx_{0}, \bfv_{0}) \in  \cK_{r}$ to be the position of the test particle at time $t$ with  initial condition $(\bfx_{0}, \bfv_{0})$.

In \cite{M-Sannals2}, the authors proved the Boltzmann-Grad limit of the periodic Lorentz gas in any dimension, which is reformulated in Theorem 3.1, \cite{M-T16} as follows. 
\bthm\label{b.g.limit}
Let $d \geq 2$. 
\benum
{\item
Assume $\bfv_{0}$ is distributed according to an absolutely continuous Borel probability measure $\lambda$ on $S^{d-1}_{1}$.  Then there exist a process $\bfQ_{n}$ on $\bbR^{d}$ with $\bfQ_{0} = 0$ such that as $r \rightarrow 0$, we have 
$$
\bfq_{n} - \bfq_{0} \Rightarrow \bfQ_{n}.
$$
}
{\item
Assume $(\bfx_{0}, \bfv_{0})$ is distributed according to an absolutely continuous Borel probability measure $\Lambda$ on $\cK_{r} \times S^{d-1}_{1}$.  Then there exist a process $\bfX_{t}$ on $\bbR^{d}$ with $\bfX_{0} = 0$ such that as $r \rightarrow 0$, we have 
$$
\bfx_{t} - \bfx_{0} \Rightarrow \bfX_{t}.
$$
}
\eenum
\ethm

Then with Theorem~\ref{superdiffusion}, the authors established the superdiffusive central limit theorem of the periodic Lorentz gas, first under the limit of $r \rightarrow 0$ and then $t \rightarrow \infty$. We stress that the order of the limits indicates that we first obtain the Boltzmann-Grad limit, and then study the long time behavior of the Boltzmann-Grad limit, which turns out to be a superdiffusion. As we have explained in the introduction, at present we do not know any results when taking the limits $r \rightarrow 0$ and $t \rightarrow \infty$ at the same time, or changing the order of the limits.

We now provide more details on the Boltzmann-Grad limit. Indeed, the displacements $\bfQ_{n}$ and $\bfX_{t}$ can be constructed from a Markov chain. To see this, define the Markov chain $(\xi_{n}, \bdeta_{n})_{n \geq 1}$ on the state space $\bbR_{> 0} \times \cB^{d}_{1}$, whose transition probability given by
\beqna\label{transition.1}
\bbP((\xi_{n}, \bdeta_{n}) \in \cA | (\xi_{n-1}, \bdeta_{n-1})) = \int_{\cA} \Phi_{0}(\bdeta_{n-1}, x, \bfz)dx d\bfz,
\eeqna
where $ \Phi_{0}(\bfomega, x, \bfz)$ is the transition kernel, which is independent of the scattering map, $\cL$ and the initial data. Here $\xi_{n} \in \bbR^{+}$ is the free path length between each collision,  and
$\bdeta_{n} \in  \cB^{d}_{1}$ is decided by the hitting position and the angle between the impacting and leaving velocities of each collision, such that the velocity after $nth$ collision $\bfV_{n}$ can be written as
\beqnas
\bfV_{n} = R(\bfv_{0})S(\bdeta_{1}) \dots  S(\bdeta_{n})\bfe_{1}, \ \ \  \bfV_{0} = \bfv_{0}; 
\eeqnas
where $R, S \in SO(d)$ are matrices decided by the scattering map. For their explicit expressions and the proof, see Lemma 2.1 in~\cite{M-T16}. 

Let 
\beqnas
\bfQ_{n}  &=&  \sum^{n}_{j=1} \xi_{j}\bfV_{j-1},\\
\bfX_{t} &=& \bfQ_{\nu_{t}} + (t - \tau_{\nu_{t}})\bfV_{\nu_{t}}, 
\eeqnas
be the distribution limits of $\bfq_{n} - \bfq_{0}$ and $\bfx_{t} - \bfx_{0}$ respectively. 

As in \cite{M-T16}, define
\beqnas
\tau_{n} = \sum^{n}_{i=1}\xi_{i}
\eeqnas
to be the time up to the nth collision. Then
\beqnas
\nu_{t} = \max\{n \in \bbZ_{\geq 0}; \tau_{n} \leq t \}
\eeqnas
is the number of collisions before time $t$. 

In \cite{M-Sannals2} it was proved that $(\xi_{1}, \bdeta_{1})$ is distributed according to 
$$
\Phi_{0}(x, \bfz) = \frac{1}{v_{d-1}} \int_{\cB^{d}_{1}}\Phi_{0}(\bfomega, x, \bfz)d\bfomega. 
$$
Moreover, by the property of $\Phi_{0}(\bfomega, x, \bfz)$, it can be verified that $\Phi_{0}(x, \bfz)$ is the the stationary measure of the process $\bfQ_{n}$. 

Since $\bfQ_{n}$ and $\bfX_{t}$ are independent of the initial data $\bfx_{0}$ and $\bfq_{0}$,  the following assumptions (noted as $\bfA$) were made on the initial data in~\cite{M-T16}, which 
\benum
{\item\label{assump1.1}
The initial velocity $\bfv_{0} \in \bbS^{d}_{1}$  is fixed,
}
{\item\label{assump1.2}
$ (\xi_{1}, \bdeta_{1})$ is distributed according to $\Phi_{0}(\bfomega, \bfz)$.
}
\eenum

Actually in \cite{M-T16} the central limit theorem can be proved for more general initial data. In this article we adopt the assumptions $\bfA$ in our main results for simplicity, although they also hold for the general initial data as in  \cite{M-T16}.

From \eqref{transition.1} we can see that $\bdeta = \{\bdeta_{n}\}^{\infty}_{n=1}$ itself is a Markov chain on $\cB^{d-1}_{1}$ with the transition probability
$$
\bbP(\bdeta_{n} \in \cA | \bdeta_{n-1}) = \int_{\cA} K_{0}(\bdeta_{n-1}, \bfz)d\bfz,
$$
where
$$
K_{0}(\bfomega, \bfz) = \int^{\infty}_{0}\Phi_{0}(\bfomega, x, \bfz)dx. 
$$
By the fact that $\int^{\infty}_{0}\int_{\cB^{d}_{1}}\Phi_{0}(\bfomega, x, \bfz)dxd\bfz = 1$ it  can be deduced that $\bdeta_{1}$ is uniformly distributed on $\cB^{d}_{1}$, and the uniform measure on $\cB^{d}_{1}$ is the  stationary measure of $\bdeta$.

Let 
$$
\Phi_{0}(x) = \frac{1}{v_{d-1}} \int_{\cB^{d}_{1}}\int_{\cB^{d}_{1}}\Phi_{0}(\bfomega, x, \bfz)d\bfomega d\bfz 
$$
be the distribution density of $\xi_{i}$. In Theorem 1.14, \cite{M-Sgafa}, the authors showed the asymptotic expansion formula of $\Phi_{0}(x)$.  See also the formulas $(5.11)$ in \cite{M-T16}. More precisely, when $d=3$ one has
\beqna\label{asym.phi1}
\Phi_{0}(x) = \Theta_{d}x^{-3} + O(x^{-3 - \frac{2}{d}}) \times \log  x,
\eeqna
and when $d = 2$, $d \geq 4$
\beqna\label{asym.phi2}
\Phi_{0}(x) = \Theta_{d}x^{-3} + O(x^{-3 - \frac{2}{d}}),
\eeqna
where $\Theta_{d} = \frac{2^{2-d}}{d(d+1)\zeta(d)}$.

%
  
\subsection{Stein's method and exchangeable pairs}
In this part, we present some elementary facts and estimates on Stein's method for multivariate normal approximation, which will be applied later in this paper. Then we give a brief description on how the  technique of exchangeable pair works. For more details, see Chen-Shao~\cite{C-S11} and Chatterjee~\cite{notes_cha}. 

Stein's method for normal approximation is a tool to estimate the difference between the distribution of a random variable $\bfw$ on $\bbR^{d}$ and the standard Gaussian random variable $\bfZ$ on $\bbR^{d}$, in the sense that
$$
\sup_{h \in \cH}|\bbE h(\bfw) - \bbE h(\bfZ)|,
$$
where $\cH$ is a given function space on $\bbR^{d}$. 

As we have explained in the introduction, the first step in Stein's method is the following Stein's equation, which turns the difference $h(\bfw) - \bbE h(\bfZ)$ into the values with respect to functions of another function space.
 More precisely, the multivariate Stein's equation for $h$ is given by
\beqna\label{stein.eq.multid}
h(\bfw) - \bbE h(\bfZ) = \Delta f_{h}(\bfw) - \bfw \cdot \nabla f_{h}(\bfw),
\eeqna
where the solution to the Stein's equation $f_{h}$ is a function on $\bbR^{d}$. Notice that the right hand side of this equation is just the O-U operator acting on $f_{h}$. 

Assume $h \in \cC^{\infty}_{0}(\bbR^{d})$, the solution to~\eqref{stein.eq.multid} can be explicitly expressed by
$$
f_{h}(\bfw) = -\int^{\infty}_{0}(\bbE(h(\bfw e^{-u} + \sqrt{1 - e^{-2u}}\bfZ)) - \bbE(h(\bfZ)))du.
$$

In the following we quote the bounds of the derivatives of $f_{h}$ from Lemma 2 in ~\cite{EMeckes} and Lemma 3 in~\cite{C-M07}.

For $k \geq 1$, $\bfx \in \bbR^{d}$,   denote the kth derivative of a function $f \in \cC^{k}(\bbR^{d})$ at $\bfx$ by $D^{k}f(\bfx)$, and the inner product $\langle D^{k}f(\bfx), (\bfu_{1}, \bfu_{2}, \dots, \bfu_{k}) \rangle$ by
$$
\langle D^{k}f(\bfx), (\bfu_{1}, \bfu_{2}, \dots, \bfu_{k}) \rangle = \sum^{d}_{i_{1}, \dots, i_{k} = 1}\frac{\partial^{k}f(\bfx)}{\partial_{x_{i_{1}}}\partial_{x_{i_{2}}}\dots\partial_{x_{i_{k}}}}(\bfu_{1})_{ i_{1}}\dots(\bfu_{k})_{ i_{k}},
$$
where $\bfu_{1}, \dots, \bfu_{k} \in \bbR^{d}$.  Define the operator norm of $D^{k}f(\bfx)$ by
$$
\|D^{k}f(\bfx)\| = \sup_{\bfu_{1}, \dots, \bfu_{k} \in \bbR^{d}}\{|\langle D^{k}f(\bfx), (\bfu_{1}, \bfu_{2}, \dots, \bfu_{k}) \rangle| , |\bfu_{1}| = \dots = |\bfu_{k}|=1\}.
$$

\blem
The following estimates hold for the solution $f_{h}$ to the Stein's equation~\eqref{stein.eq.multid} with respect to $h$. 
\bitem 
{\item
If $h \in \cC^{1}(\bbR^{d})$, we have
\beqna\label{bound.deri1}
\sup_{\bfx \in \bbR^{d}}\|Df_{h}(\bfx)\| \leq \sqrt{\frac{\pi}{2}}\sup_{\bfx \in \bbR^{d}}\|Dh(\bfx)\|,
\eeqna
}
{\item
If $h \in \cC^{1}(\bbR^{d})$, we have
\beqna\label{bound.deri2}
\sup_{\bfx \in \bbR^{d}}\|{\rm Hess}f_{h}(\bfx)\|_{H.S.}  \leq  \sup_{\bfx \in \bbR^{d}}\|Dh(\bfx)\|,
\eeqna
where for a matrix $A$, $\|A\|_{H.S.} := \sqrt{{\rm Tr}(AA^{t})}$.
}
{\item
If  $h \in \cC^{2}(\bbR^{d})$, we have
\beqna\label{bound.deri3}
\sup_{\bfx \in \bbR^{d}}\|D^{3}f_{h}(\bfx)\| \leq  \frac{\sqrt{2\pi}}{4}\sup_{\bfx \in \bbR^{d}}\|D^{2}h(\bfx)\|.
\eeqna
}
\eitem
\elem

Now we explain how the method of exchangeable pairs works. We start with the definition of the exchangeable pairs. 

\bdefi[exchangeable pairs]\label{def.exc.pair}
Let $\bfW$ and $\bfW'$ be two identically distributed random vectors on $\bbR^{d}$ on the same probability space. We call $(\bfW, \bfW')$ an exchangeable pair if $(\bfW, \bfW')$ has the same distribution as $(\bfW', \bfW)$. 
\edefi

In the univariate case, to apply the method of exchangeable pairs, we need the exchangeable pair $(W, W')$ to satisfy the "linear regression condition", i.e.
\beqna\label{ep.linear}
\bbE(W' - W| W) = -\lambda W
\eeqna
holds for some constant $\lambda$. In the multivariate case, there are various analogues to this linear condition~\eqref{ep.linear}, see for example \cite{C-M07} and \cite{RR09}.  In this paper we employ the conditions in Theorem 4, \cite{C-M07}, which is enough for our case. More precisely, there exist  constants $\lambda$, $\sigma$ and a random matrix $E$, such that the following two conditions hold
\bitem
{\item
\beqna\label{exc.pair.eq1}
\bbE(\bfW' - \bfW| \bfW) = -\lambda \bfW,
\eeqna
}
{\item
\beqna\label{exc.pair.eq2}
\bbE((\bfW' - \bfW)(\bfW' - \bfW)^{t}| \bfW) = 2\lambda \sigma^{2}\Id_{d} + 2\lambda\bbE(E|\bfW).
\eeqna
}
\eitem

Notice that by definition,  the exchangeable pair $(\bfW, \bfW')$ has the same distribution as that of $(\bfW', \bfW)$. A direct consequence is that for any anti-symmetric function $g$ on $\bbR^{d} \times \bbR^{d}$, we have
$$
\bbE g(\bfW, \bfW') = 0. 
$$

To explain how the exchangeable pair works, we do the following heuristic computation. Let $f$ be a function on $\bbR^{d} \times \bbR^{d}$, given by
\beqnas
f(\bfx, \bfx') = \langle (\bfx' -  \bfx), \nabla f_{h}(\bfx') +  \nabla f_{h}(\bfx) \rangle,
\eeqnas
where $f_{h}$ is the solution to the Stein's equation~\eqref{stein.eq.multid} and we assume that $f_{h} \in \cC^{3}_{b}$. Since $f$ is anti-symmetric, we derive that
\beqna\label{eq2.2.2}
\nonumber 0 &=& \bbE [\langle (\bfW' -  \bfW), \nabla f_{h}(\bfW') +  \nabla f_{h}(\bfW) \rangle]\\
\nonumber &=& \bbE [\langle (\bfW' -  \bfW), \nabla f_{h}(\bfW') -  \nabla f_{h}(\bfW) \rangle + 2\langle (\bfW' -  \bfW),  \nabla f_{h}(\bfW) \rangle]\\
&=& \bbE [\langle (\bfW' -  \bfW)(\bfW' -  \bfW)^{t}, {\rm Hess}f_{h}(\bfW) \rangle + 2\langle (\bfW' -  \bfW),  \nabla f_{h}(\bfW) \rangle + R], 
\eeqna
where $R$ is the error in Taylor expansion, given by
\beqnas
R  &=&  \frac{1}{2} \sum^{d}_{j, k, l=1}D^{3}_{jkl}f_{h}(\bfW + \tau(\bfW' -  \bfW))(\bfW' -  \bfW)_{j}(\bfW' -  \bfW)_{k}(\bfW' -  \bfW)_{l},
\eeqnas
for some $\tau \in (0, 1)$.
Thus by Stein's equation~\eqref{stein.eq.multid}, we have
\beqnas
& &\bbE h(\bfW) - \bbE h(\bfZ) \\
&=& \bbE [\Delta f_{h}(\bfW) - \bfW \cdot \nabla f_{h}(\bfW)]\\
&=&  \bbE [\langle \Id - \bbE((\bfW' -  \bfW)(\bfW' -  \bfW)^{t}|\bfW), {\rm Hess}f_{h}(\bfW) \rangle]  \\
& &+ \bbE [\langle (\bfW' -  \bfW)(\bfW' -  \bfW)^{t}, {\rm Hess}f_{h}(\bfW) \rangle]\\
& &  + \bbE[\langle -2\bbE((\bfW' -  \bfW)|\bfW) - \bfW, \nabla f_{h}(\bfW) \rangle + 2\langle (\bfW' -  \bfW),  \nabla f_{h}(\bfW) \rangle ]\\
&=&  \bbE [\langle \Id - \bbE((\bfW' -  \bfW)(\bfW' -  \bfW)^{t}|\bfW), {\rm Hess}f_{h}(\bfW) \rangle ]\\
& & + \bbE[ \langle -2\bbE((\bfW' -  \bfW)|\bfW) - \bfW, \nabla f_{h}(\bfW) \rangle] - \bbE R\\
&=&  \bbE [\langle (1- 2\lambda \sigma^{2})\Id_{d} - 2\lambda\bbE(E|\bfW), {\rm Hess}f_{h}(\bfW) \rangle ]\\
& & + \bbE[ \langle (2\lambda  - 1)\bfW, \nabla f_{h}(\bfW) \rangle] - \bbE R,
\eeqnas
where in the last line is due to \eqref{eq2.2.2} and  the conditions \eqref{exc.pair.eq1}, \eqref{exc.pair.eq2}. Thus we conclude that the estimate of $|\bbE h(\bfW) - \bbE h(\bfZ)|$ boils down to obtain the bounds on $M_{1}f_{h}$, $M_{2}f_{h}$, $M_{3}f_{h}$ and $\bbE(E|\bfW)$, and also the moment estimates of $\bfW$. 

We point out that the above discussions could not directly apply to the Boltzmann-Grad limit of the periodic Lorentz gas. In the proof of the main Theorem~\ref{thm1} (Section 4), we conduct the truncation procedure and then take the conditional expectations before employing the Stein's method.

\section{Key estimates}
In this section we provide the key estimates to control the error term in Stein's method.  They are analogous to the moment bounds in the classical settings, as we explain in Section 2.2. 

First we adopt the truncation procedure in \cite{M-T16}. More precisely, for $1 \leq i \leq n$, let  
$r^{2}_{i, n}=   n(\log n)^{\gamma}$ with $0 < \gamma < 1$ and define
$$
\xi'_{i,n} = \xi_{i}\un_{\{\xi_{i}^{2} \leq r^{2}_{i, n}\}}
$$
and 
$$
\bfQ'_{n} = \sum^{n}_{i=1}\xi'_{i, n}\bfV_{i-1}, \ \ \ \bfW'_{n} = \frac{1}{\sqrt{n \log n}}\bfQ'_{n}.
$$ 

Notice that by the asymptotic expansion formula~\eqref{asym.phi1} and  ~\eqref{asym.phi2} of $\Phi_{0}$ we deduce that
\beqna\label{estimate.qn}
\nonumber& &\sup_{n}\bbE \|\bfQ_{n} - \bfQ'_{n}\| \leq \bbE(\xi_{1})+\sum^{n}_{i=2} \bbE(\xi_{i} \un_{\{\xi^{2}_{i} \geq r^{2}_{i, n}\}}) \\
\nonumber & \leq&  1 + \sum^{n}_{i=2} \int^{\infty}_{n^{\frac{1}{2}}(\log n)^{\frac{\gamma}{2}}} (x^{-2} + O(x^{-2 - \frac{2}{d}}\log x))dx\\
 &=&  O( \frac{n^{\frac{1}{2}}}{(\log n)^{\frac{\gamma}{2}}}),
\eeqna
so that
\beqna\label{estimate.wn}
\sup_{n}\bbE \|\bfW_{n} - \bfW'_{n}\| = O( \frac{1}{(\log n)^{\frac{\gamma + 1}{2}}}) .
\eeqna

As in \cite{M-T16}, we divide $\xi'_{i, n}$ into two parts, 
$$
\xi'_{i, n} = \tilde{\xi}_{i, n} + m_{i, n}, 
$$
where
$$
m_{i, n} = \bbE(\xi'_{i, n}|\bdeta),  \ \ \ \ \   \tilde{\xi}_{i, n} = \xi'_{i, n} - m_{i, n}.
$$

The rest of this section is devoted to the moment estimates and the key estimates, which are due to the spectral gap of the transition kernels. 

\subsection{The moment estimates}
In this section, we prove the moment estimates for $\xi'_{i, n}$, $m_{i,n}$ and $\bbE(\xi'^{k}_{i, n}|\bdeta)$. These estimates are the consequence of the asymptotic behavior of $\Phi_{0}(\bfomega, x, \bfz)$, which was proved in \cite{M-Sgafa}.  Although some of the estimates have appeared in \cite{M-T16} ,  for the completeness of this paper we quote them here without proof. Notice that our truncation bound $r_{i, n}$ is slightly different from theirs, but the proof is the same.

%

\blem
\beqna\label{m.est.xin.2}
\bbE \left( (\xi_{i, n}')^{2}\right) =d\sigma^{2}_{d}\log n + O(\log \log n),
\eeqna
where the constant $\sigma^{2}_{d} = \frac{2^{1-d}}{d^{2}(d+1)\zeta(d)}$, and $\zeta(d) = \sum^{d}_{n=1}n^{-d}$ is the Riemann zeta function.

For $k \geq 3$, we have
\beqna\label{m.est.xin.k}
\bbE \left( (\xi_{i, n}')^{k} \right) = O(r^{k-2}_{i, n}) .
\eeqna

\elem
\bpf
This is a direct consequence of the asymptotic expansion formulas~\eqref{asym.phi1} and ~\eqref{asym.phi2} of $\Phi_{0}(x)$. Indeed, we have
\beqnas
\bbE(({\xi}^{'}_{i})^{k}) &=& \int^{r_{i, n}}_{0} x^{k}\Phi_{0}(x)dx =  \int^{r_{i, n}}_{0} x^{k}(\Theta_{d}x^{-3} + O(x^{-3 - \frac{2}{d}}))dx\\
&=& O(r^{k-2}_{i, n}).
\eeqnas

Also see Proposition 6.10 in \cite{M-T16} for \eqref{m.est.xin.2}.  

\epf


\blem
For $k \geq 1$, we have
\beqnas\label{m.est.m.k}
|\bbE(m^{k}_{i, n})| < O(n^{\frac{k-2}{2}}(\log n)^{\frac{\gamma(k-2)}{2}}).
\eeqnas
In particular, we have 
\beqna\label{m.est.m.1}
\bbE m_{i, n} = \bar{\xi} + O(\frac{1}{r_{i, n}}),
\eeqna
where $\bar{\xi} = \frac{\Gamma(\frac{d+1}{2})}{\pi^{\frac{d-1}{2}}}$ is the limit of the mean free path length as $r \rightarrow 0$.

As for $d=2$, we have
\beqna\label{m.est.m.2.d2}
\bbE m^{2}_{i, n} = O(\log \log n),
\eeqna
while for $d \geq 3$, 
\beqna\label{m.est.m.2.d3}
\bbE m^{2}_{i, n} = O(1).
\eeqna
\elem

\bpf
The estimates \eqref{m.est.m.1}, \eqref{m.est.m.2.d2} and \eqref{m.est.m.2.d3} are proven in Proposition 6.8, 6.9 in \cite{M-T16}. 

Recall Prop 6.7 in \cite{M-T16} that
\beqnas
\bbP(m_{i, n} > u) = O(u^{-(1+\frac{d}{2})})\un_{u \leq r_{i, n}}. 
\eeqnas
Thus we derive that
\beqnas
|\bbE(m^{k}_{i, n})| &=& |\int^{r_{i, n}}_{0} ku^{k-1}\bbP(m_{i} > u)du|
 = O(r_{i, n}^{k-1-\frac{d}{2}}) < O(n^{\frac{k-2}{2}}(\log n)^{\frac{\gamma(k-2)}{2}}).
\eeqnas

\epf

\blem
\beqna
\label{m.est.tilde.2.2}\bbE(\bbE(\tilde{\xi}_{i, n}^{2}|\bdeta)^{2}) &\leq&  O(r^{2}_{i, n}) = O(n(\log n)^{\gamma}),\\
\label{m.est.tilde.3.2}\bbE(\bbE(\tilde{\xi}_{i, n}^{3}|\bdeta)^{2})  &\leq&    O(r^{4}_{i, n}) = O(n^{2}(\log n)^{2\gamma}).
\eeqna

%
\elem
\bpf
%

The first estimate~\eqref{m.est.tilde.2.2} is proved in Lemma 6.11 in \cite{M-T16}. 

As for the second one, notice that
\beqna\label{l3.3.0}
\bbE(\bbE(\tilde{\xi}_{i, n}^{3}|\bdeta)^{2}) &=& \bbE(\bbE((\xi'_{i, n} - m_{i, n})^{3}|\bdeta)^{2}) \\
\nonumber&=& \bbE \big( \bbE((\xi'_{i, n})^{3}|\bdeta) - 3\bbE((\xi'_{i, n})^{2}|\bdeta)m_{i, n} + 2m^{3}_{i,n}\big)^{2}\\
\nonumber&\leq& C\bbE \big( \bbE((\xi'_{i, n})^{3}|\bdeta)^{2} +\bbE((\xi'_{i, n})^{2}|\bdeta)^{2}m^{2}_{i, n} + m^{6}_{i, n}\big)\\
\nonumber&\leq& C\big(\bbE(\bbE((\xi'_{i, n})^{3}|\bdeta)^{2}) + \bbE(\bbE((\xi'_{i, n})^{2}|\bdeta)^{3})  + \bbE m^{6}_{i, n}\big),
\eeqna
where the last line is due to Young's inequality, i.e.
$$
\bbE (\bbE((\xi'_{i, n})^{2}|\bdeta)^{2}m^{2}_{i, n}) \leq   \frac{2}{3}\bbE (\bbE((\xi'_{i, n})^{2}|\bdeta)^{3} + \frac{1}{3}\bbE m^{6}_{i, n}.
$$


Thanks to ~\eqref{m.est.m.k}, we obtain
\beqna\label{l3.3.1}
\bbE m^{6}_{i, n} \leq O(r^{4}_{i, n}) = O(n^{2}(\log n)^{2\gamma}).
\eeqna

Now we turn to the estimates of $\bbE(\bbE((\xi'_{i})^{3}|\bdeta)^{2})$ and $\bbE(\bbE((\xi'_{i})^{2}|\bdeta)^{3})$. 

For $\bbE(\bbE((\xi'_{i})^{3}|\bdeta)^{2})$, we have 
\beqna\label{l3.3.2}
\nonumber \bbE(\bbE((\xi'_{i})^{3}| \bdeta)^{2})  
&=& \frac{1}{v_{d-1}}\int_{\cB^{d-1}_{1}}\int_{\cB^{d-1}_{1}}(\int^{r_{i}}_{0}x^{3}\Phi_{0}(\bfomega, x, \bfz)dx)^{2}\frac{1}{K_{0}(\bfomega, \bfz)}d\bfomega d\bfz\\
\nonumber &\leq&  C\int_{\cB^{d-1}_{1}}\int_{\cB^{d-1}_{1}}(\int^{r_{i}}_{0}x^{6}\Phi_{0}(\bfomega, x, \bfz)dx)(\int^{r_{i}}_{0}\Phi_{0}(\bfomega, x, \bfz)dx)\frac{1}{K_{0}(\bfomega, \bfz)}d\bfomega d\bfz\\
&=&  O(r^{4}_{i, n}) = O(n^{2}(\log n)^{2\gamma}).
\eeqna

Then $\bbE (\bbE((\xi'_{i})^{2}|\bdeta)^{3})$ can be estimated in the same way. 
\beqna\label{l3.3.3}
\bbE (\bbE((\xi'_{i})^{2}|\bdeta)^{3})
&=&  \frac{1}{v_{d-1}}\int_{\cB^{d-1}_{1}}\int_{\cB^{d-1}_{1}}(\int^{r_{i}}_{0}x^{2}\Phi_{0}(\bfomega, x, \bfz)dx)^{3}\frac{1}{K_{0}(\bfomega, \bfz)^{2}}d\bfomega d\bfz\\
\nonumber&\leq&  O(r^{4}_{i, n}) = O(n^{2}(\log n)^{2\gamma}).
\eeqna
Combining \eqref{l3.3.1}, \eqref{l3.3.2} and \eqref{l3.3.3}, we deduce \eqref{m.est.tilde.3.2}.


%

\epf

\blem
\beqna\label{m.est.tilde.2}
\bbE \tilde{\xi}_{i, n}^{2}  
= d\sigma^{2}_{d}\log n +  O(\log \log n), 
\eeqna
and 
\beqna\label{m.est.tilde.4}
\bbE \tilde{\xi}_{i, n}^{4}  = O(r^{2}_{i, n}) = O(n(\log n)^{\gamma}).
\eeqna
\elem

\bpf
Notice that
\beqnas
\bbE \tilde{\xi}_{i, n}^{4}  = \bbE (\xi^{'}_{i, n} - m_{i, n})^{4} \leq  C(\bbE (\xi^{'}_{i, n})^{4} + \bbE(m_{i, n}^{4})), 
\eeqnas
such that by \eqref{m.est.xin.k} and \eqref{m.est.m.k}, we obtain \eqref{m.est.tilde.4}.
\epf

\subsection{The key estimates}
In this section our aim is to estimate expectations in the form of
$$
\bbE(\sum^{n}_{i=1}(e_{j} \cdot \bfV_{i})(e_{k} \cdot \bfV_{i})(e_{l} \cdot \bfV_{i})f(\bfV_{i}, \dots,\bfV_{i+m})),
$$ 
which will appear in the error term when applying Stein's method.  Here $1 \leq j, k, l \leq d$, $m \in \bbN$. We will define the precise function space of $f$ below. 

Recall that in \cite{M-T16}, the authors proved an exponential mixing estimates and gave estimates of $\bbE((e_{i} \cdot \bfV_{n})f(\bdeta_{n}, \dots, \bdeta_{n +m}))$ and $\bbE((e_{i} \cdot \bfV_{n})^{2}f(\bdeta_{n}, \dots, \bdeta_{n +m}))$, as a consequence of the spectral gap of the transition operator of $\bdeta$. More precisely, they define the operator $P$ on the Hilbert space $\cH = L^{2}(\cB^{d-1}_{1},V, \frac{1}{v_{d-1}}d\bfomega)$,
$$
Pf(\bfomega) = \bbE(f(\bdeta_{n})|\bdeta_{n-1} = \bfomega) = \int_{\cB^{d-1}_{1}} K_{0}(\bfomega, \bfz)f(\bfz)d\bfz.
$$ 
where $V$ is a finite dimension vector space, and prove that $P$ has the spectral gap $1 - \omega_{0}$ with $\omega_{0} < 1$ by Doeblin theory.

In this paper we consider $\bdalpha_{n} = (\bfV_{n-1}, \bfV_{n}, \bfV_{n+1})$ as a stationary Markov chain, which is a direct consequence of the results in~\cite{M-Sannals2}, as pointed out in \cite{M-T16}. It was proved in Prop 13.4, \cite{M-T16} that $\bdalpha_{n}$ also admits a spectral gap. Indeed, define
$$
\cV = \{(\bfv_{n-1}, \bfv_{n}, \bfv_{n+1}) \in (S^{d-1}_{1})^{3}, \phi(\bfv_{n-1}, \bfv_{n}) >  B_{\theta}, \phi(\bfv_{n}, \bfv_{n+1}) >  B_{\theta}\}, 
$$
where $\phi(\bfu_{1}, \bfu_{2}) \in [0, \pi]$ denotes the angle between two vectors $\bfu_{1}, \bfu_{2}$, and $B_{\theta} = \inf_{\omega \in [0, 1)}|\theta(\omega)|$ is the minimal deflection angle associated with the scattering angle $\theta(\omega)$. To ensure that the Boltzmann-Grad limit of the periodic Lorentz gas exists, see assumptions $(A)$, $(B)$ on $\theta$ in Section 2, \cite{M-T16}.

Let  $\cP_{\bfv}$ be the  transition operator of $\bdalpha_{n}$, given by
\beqnas
\cP_{\bfv} f(\bfz) = \bbE (f(\bdalpha_{n})|\bdalpha_{n-1} = \bfz) = \int_{\cV} \cK_{\bfv}(\bfz, \bfomega)f(\bfomega)d\mu(\bfomega),
\eeqnas
where $f \in L^{2}(\cV, d\mu)$, $\cK_{\bfv}$ is the transition kernel of $\bdalpha_{n}$ and $d\mu(\bdalpha_{n})$ is its stationary measure. For the explicit expressions of $d\mu(\bfomega)$ and $\cK_{\bfv}(\bfz, \bfomega)$, see \cite{M-T16} and \cite{M-Sannals2}. 

It should be pointed out that both $\cK_{\bfv}(\bfz, \bfomega)$ and $d\mu$ are invariant under rotation, i.e. for $R \in SO(d)$, we have
$$
\cK_{\bfv}(R\bfz, R\bfomega) = \cK_{\bfv}(\bfz, \bfomega), \ \ \  d\mu(R\bfomega) = d\mu(\bfomega).
$$


Moreover,  it is shown in \cite{M-T16} that the classical Doeblin theory leads to the spectral gap of $\cP_{\bfv}$, denoted by $1 - \omega_{\bfv}$ with $\omega_{\bfv} = \|\cP_{\bfv} - \Pi_{\bfv}\| < 1$, where $\|\cdot \|$ is the operator norm in $L^{2}(\cV, \mu)$. In addition, define
\beqnas
\Pi_{\bfv} f(\bfomega) = \int_{\cV}f(\bfz)d\mu(\bfz),
\eeqnas
then we have
\beqna\label{exchangable}
\Pi_{\bfv} \cP_{\bfv} = \cP_{\bfv} \Pi_{\bfv} = \Pi_{\bfv}.
\eeqna
%
%
%


Let $L^{2}_{0}(\cV, \mu)$ be the orthogonal complement of the constant functions in $L^{2}(\cV, \mu)$. Then for  any $f \in L^{2}_{0}(\cV, \mu)$, there exists $g \in L^{2}(\cV, \mu)$ such that   $f = (\cI - \cP_{\bfv})g$.

Denote by $\|f\|^{2}_{2} = \int_{\cV}|f(\bfz)|^{2}d\mu(\bfz)$ the norm of $L^{2}(\cV, \mu)$. The following lemma is a direct consequence of spectral gap, of which we omit the proof.

\blem\label{keylemma}
For any $f \in L^{2}(\cV, \mu)$,  we have
\beqna\label{keylemma.ineq}
\|\cP^{n}_{\bfv}f - \Pi_{\bfv} f\|_{2} \leq C\omega^{n}_{\bfv}\|f\|_{2}.
\eeqna
\elem


Now we turn to the main estimate of this section. We say a function  $f: \cV^{m+1} \rightarrow \bbR$ is rotation-invariant if for any $R \in SO(d)$, 
$$
f(R\bdalpha_{0}, \dots, R\bdalpha_{m}) = f(\bdalpha_{0}, \dots, \bdalpha_{m}),
$$
where $R\bdalpha_{n} := (R\bfV_{n-1}, R\bfV_{n}, R\bfV_{n+1})$.

\bprop\label{triple-estimate}
For measurable and rotation-invariant functions  $f: \cV^{m+1} \rightarrow \bbR$ with 
$$
\|f\|_{2}^{2} := \bbE(f(\bdalpha_{0}, \dots, \bdalpha_{m})^{2}) < \infty,
$$
we have
\beqna\label{estimate.deux2}
\sum^{n}_{i=1}|\bbE(\big((\bfe_{j} \cdot \bfV_{i})^{2} - \frac{1}{d}\big)f(\bdalpha_{i}, \dots, \bdalpha_{i +m}))| \leq C_{2, 2}\|f\|_{2}, 
\eeqna 
where $C_{2}$ is a constant independent of $n$, and for $j \neq k$, 
\beqna\label{estimate.deux1}
\bbE((\bfe_{j} \cdot \bfV_{i}) (\bfe_{k} \cdot \bfV_{i})f(\bdalpha_{i}, \dots, \bdalpha_{i +m})) = 0.
\eeqna

Moreover, we have
\beqna\label{estimate.trois}
\bbE ((\bfe_{j} \cdot \bfV_{i})(\bfe_{k} \cdot \bfV_{i})(\bfe_{l} \cdot \bfV_{i})f(\bdalpha_{i}, \dots, \bdalpha_{i +m}))  = 0.
\eeqna
\eprop

\bpf
We begin with the simplest case. Notice that
\beqnas
\bbE((\bfe_{j} \cdot \bfV_{i})f(\bdalpha_{i}, \dots, \bdalpha_{i +m})) = \bbE ((\bfe_{j} \cdot \bfV_{i})\bbE(f(\bdalpha_{i}, \dots, \bdalpha_{i +m})|\bdalpha_{i})).
\eeqnas
Let $ \tilde{f}(\bdalpha_{i}) := \bbE(f(\bdalpha_{i}, \dots, \bdalpha_{i +m})|\bdalpha_{i}))$, then by the rotation invariances of $f$, $\cK(\bfz, \bfomega)$ and $d\mu$ we conclude that $\tilde{f}$ is also rotation-invariant, i.e. for any $R \in SO(d)$, 
\beqnas
\tilde{f}(\bdalpha_{i}) = \tilde{f}(R\bdalpha_{i}). 
\eeqnas
By taking $R = -\Id$, it is easy to verify that
\beqna\label{estimate.une}
\bbE((\bfe_{j} \cdot \bfV_{i})f(\bdalpha_{i}, \dots, \bdalpha_{i +m})) = \int_{\cV} (\bfe_{j} \cdot \bfV_{i})\tilde{f}(\bdalpha_{i})d\mu(\bdalpha_{i})= 0.
\eeqna

Now we prove the estimates~\eqref{estimate.deux1}, \eqref{estimate.deux2}. By the rotation-invariance, we have for any $R \in SO(d)$, 
\beqnas
\int_{\cV} (\bfe_{j} \cdot \bfV_{i})^{2}\tilde{f}(\bdalpha_{i})d\mu(\bdalpha_{i}) 
= \int_{\cV} (\bfe_{j} \cdot R\bfV_{i})^{2}\tilde{f}(\bdalpha_{i})d\mu(\bdalpha_{i}) = \int_{\cV} (R^{t}\bfe_{j} \cdot \bfV_{i})^{2}\tilde{f}(\bdalpha_{i})d\mu(\bdalpha_{i}),
\eeqnas
thus by choosing $R^{t}_{k} \in SO(d)$ for $k = 1,\dots, d$, such that $R^{t}_{k}\bfe_{j} = \bfe_{k}$ and the fact $\|\bfV_{i}\| =1$, we derive that
\beqna\label{eq3.6.1}
\nonumber & &\int_{\cV} (\bfe_{j} \cdot \bfV_{i})^{2}\tilde{f}(\bdalpha_{i})d\mu(\bdalpha_{i}) \\
 &=& \frac{1}{d}\sum^{d}_{k=1}\int_{\cV} (\bfe_{k} \cdot \bfV_{i})^{2}\tilde{f}(\bdalpha_{i})d\mu(\bdalpha_{i}) = \frac{1}{d}\int_{\cV} \tilde{f}(\bdalpha_{i})d\mu(\bdalpha_{i}).
\eeqna
Thus let $g(\bdalpha_{i}) = (\bfe_{j} \cdot \bfV_{i})^{2}\tilde{f}(\bdalpha_{i})$ and with \eqref{eq3.6.1} we obtain
\beqnas
& &\bbE(\big((\bfe_{j} \cdot \bfV_{i})^{2} - \frac{1}{d}\big)f(\bdalpha_{i}, \dots, \bdalpha_{i +m}))\\
&=&\bbE(\bbE((\bfe_{j} \cdot \bfV_{i})^{2}\tilde{f}(\bdalpha_{i})|\bdalpha_{0})) - \frac{1}{d}\bbE f(\bdalpha_{i}, \dots, \bdalpha_{i +m})\\
&=&\bbE\big(\cP^{i}_{\bfv}g(\bdalpha_{0}) -   \Pi_{\bfv}g \big).
\eeqnas

Therefore, Lemma~\ref{keylemma} directly leads to
\beqnas
& &|\bbE(\big((\bfe_{j} \cdot \bfV_{i})^{2} - \frac{1}{d}\big)f(\bdalpha_{i}, \dots, \bdalpha_{i +m}))|\\
&=& |\bbE\big(\cP^{i}_{\bfv}g(\bdalpha_{0}) -   \Pi_{\bfv}g \big)| \leq  \sqrt{ \bbE|\cP^{i}_{\bfv}g(\bdalpha_{0}) -   \Pi_{\bfv}g |^{2}}\\
&\leq& C\omega^{i}_{\bfv}\|g\|_{2} \leq   C\omega^{i}_{\bfv}\|f\|_{2},
\eeqnas
implying  that
\beqnas
\sum^{n}_{i=1}|\bbE(\big((\bfe_{j} \cdot \bfV_{i})^{2} - \frac{1}{d}\big)f(\bdalpha_{i}, \dots, \bdalpha_{i +m}))| \leq C_{2}\|f\|_{2}.
\eeqnas

The formula~\eqref{estimate.deux2} can be obtained as soon as we notice that 
 for $j \neq k$ the rotation invariance gives rise to
\beqnas
\int_{\cV} (\bfe_{j} \cdot \bfV_{i})(\bfe_{j} \cdot \bfV_{i})\tilde{f}(\bdalpha_{i})d\mu(\bdalpha_{i}) 
&=& \int_{\cV} (\bfe_{j} \cdot R\bfV_{i})(\bfe_{j} \cdot R\bfV_{i})\tilde{f}(\bdalpha_{i})d\mu(\bdalpha_{i}).
\eeqnas
Then by choosing $R \in SO(d)$ such that 
\beqnas
R^{t} \bfe_{i} = \bfe_{i}, \ \ \ R^{t} \bfe_{j} = -\bfe_{j},
\eeqnas
we get 
$$
\int_{\cV} (\bfe_{i} \cdot \bfV_{i})(\bfe_{j} \cdot \bfV_{i})\tilde{f}(\bdalpha_{i})d\mu(\bdalpha_{i})  = -\int_{\cV} (\bfe_{i} \cdot \bfV_{i})(\bfe_{j} \cdot \bfV_{i})\tilde{f}(\bdalpha_{i})d\mu(\bdalpha_{i}),
$$
leading to
$$
\bbE ((\bfe_{i} \cdot \bfV_{i})(\bfe_{j} \cdot \bfV_{i})f(\bdalpha_{i}, \dots, \bdalpha_{i +m})) = 0.
$$

As for~\eqref{estimate.trois}, notice that the rotation invariance leads to
\beqnas
\int_{\cV} (\bfe_{j} \cdot \bfV_{i})(\bfe_{k} \cdot \bfV_{i})(\bfe_{l} \cdot \bfV_{i})\tilde{f}(\bdalpha_{i})d\mu(\bdalpha_{i}) = \int_{\cV} (\bfe_{j} \cdot R\bfV_{i})(\bfe_{k} \cdot R\bfV_{i})(\bfe_{l} \cdot R\bfV_{i})\tilde{f}(\bdalpha_{i})d\mu(\bdalpha_{i}),
\eeqnas
such that by choosing $R$ to satisfy
\beqnas
R^{t} \bfe_{i} = \lambda_{i}\bfe_{i}, \ \ \ R^{t} \bfe_{j} = \lambda_{j}\bfe_{j}, \ \ \  R^{t} \bfe_{k} = \lambda_{k}\bfe_{k},
\eeqnas
and $\lambda_{i}\lambda_{j}\lambda_{k} = -1$, we finish the proof of \eqref{estimate.trois}.

\epf

\brmq
We point out that the estimate~\eqref{estimate.deux2} is proved in Prop 8.2,  \cite{M-T16} as a corollary of the exponential mixing estimate, and the special form of the equality~\eqref{estimate.une} is proved in display $(13.46)$ in \cite{M-T16}.  The exponential mixing estimate can also be proved by our arguments. We now claim the following argument, which will be used later in this paper: for measurable functions $f, g: \cB^{d-1} \rightarrow \bbR$ satisfying
$$
\|f\|_{2} \leq \infty, \ \ \ \ \ \|g\|_{2} \leq \infty, 
$$
we have for some constant  $\omega \in (0, 1)$
\beqna\label{est.expo.mixing}
|\Cov((\bfe, \bfV_{n_{1}})f(\bdalpha_{n_{1}}), (\bfe, \bfV_{n_{2}})g(\bdalpha_{n_{2}}))| \leq C\omega^{|n_{1} - n_{2}|}\|f\|_{2}\|g\|_{2}, 
\eeqna
where $n_{1}, n_{2} \in \bbN$.

To this end, it suffices to assume that $n_{1} \leq n_{2}$. It should be pointed out that since $\bdeta_{i}$ is the function of the angles between $\bfV_{i-1}$, $\bfV_{i}$, and $\bfV_{i}$, $\bfV_{i+1}$, it is a rotation-invariant function of $\bdalpha_{i}$.

By~\eqref{estimate.une} we see that
$$
\bbE((\bfe, \bfV_{n_{1}})f(\bdalpha_{n_{1}})) = \bbE((\bfe, \bfV_{n_{2}})g(\bdalpha_{n_{2}})) = 0.
$$
Thus we have
\beqnas
& &\Cov((\bfe, \bfV_{n_{1}})f(\bdeta_{n_{1}}), (\bfe, \bfV_{n_{2}})g(\bdeta_{n_{2}})) = \bbE((\bfe, \bfV_{n_{1}})f(\bdeta_{n_{1}})(\bfe, \bfV_{n_{2}})g(\bdeta_{n_{2}}))\\
&=& \bbE((\bfe, \bfV_{n_{1}})f(\bdeta_{n_{1}})\bbE((\bfe, \bfV_{n_{2}})g(\bdeta_{n_{2}})|\bdeta_{n_{1}}))\\
&=& \bbE((\bfe, \bfV_{n_{1}})f(\bdeta_{n_{1}})\bbE((\bfe, \bfV_{n_{2} - n_{1}+1})g(\bdeta_{n_{2} - n_{1}+1})|\bdeta_{1})),
\eeqnas
where the last line is due to the stationarity of the $\{\xi_{i}, \bdeta_{i}\}$. Recall that $\bdeta_{1}$ is uniformly distributed on $\cB^{d-1}_{1}$. Then by Lemma~\ref{keylemma} and Cauchy-Schwarz inequality we have
\beqnas
& &|\Cov((\bfe, \bfV_{n_{1}})f(\bdeta_{n_{1}}), (\bfe, \bfV_{n_{2}})g(\bdeta_{n_{2}}))| \leq   \|f(\bdeta_{n_{1}})\|_{2}\|\bbE((\bfe, \bfV_{n_{2} - n_{1}})g(\bdeta_{n_{2} - n_{1}})|\bdeta_{1})\|_{2}\\
&\leq& C\omega^{|n_{1} - n_{2}|}\|f\|_{2}\|g\|_{2}.
\eeqnas
\ermq

%
%
%

As applications of the exponential mixing estimate~\eqref{est.expo.mixing}, we prove the following estimates which we will need to obtain the convergence rate.

\blem
When $d=2$, we have
\beqna\label{est.sum.min.d2}
\|\sum^{n}_{i=1}m_{i, n}\bfV_{i-1}\|_{2} = O(\sqrt{n \log \log n}),
\eeqna
while for $d \geq 3$, we get
\beqna\label{est.sum.min.d3}
\|\sum^{n}_{i=1}m_{i, n}\bfV_{i-1}\|_{2} = O(\sqrt{n}).
\eeqna
\elem

\bpf
It suffices to give the estimate of $\|\sum^{n}_{i=1}m_{i, n}\langle\bfV_{i-1}, \bfe_{k}\rangle\|_{2}$ for any $1 \leq k \leq d$.
\beqnas
& &\sqrt{\bbE(\sum^{n}_{i=1}m_{i, n}\langle\bfV_{i-1}, \bfe_{k}\rangle)^{2}}\\
&\leq& \sqrt{\bbE(\sum^{n}_{i=1}m^{2}_{i, n}\langle\bfV_{i-1}, \bfe_{k}\rangle^{2})} + \sqrt{\bbE(\sum^{n}_{i \neq j}m_{i, n}m_{j, n}\langle\bfV_{i-1}, \bfe_{k}\rangle \langle\bfV_{j-1}, \bfe_{k}\rangle)}.
\eeqnas
For the first part, we have
$$
\bbE(\sum^{n}_{i=1}m^{2}_{i, n}\langle\bfV_{i-1}, \bfe_{k}\rangle^{2}) \leq \sum^{n}_{i=1}\bbE m^{2}_{i, n},
$$
then when $d =2$, by~\eqref{m.est.m.2.d2} we get
$$
\bbE(\sum^{n}_{i=1}m^{2}_{i, n}\langle\bfV_{i-1}, \bfe_{k}\rangle^{2}) = O(n \log \log n), 
$$
and when $d \geq 3$, \eqref{m.est.m.2.d3} leads to
$$
\bbE(\sum^{n}_{i=1}m^{2}_{i, n}\langle\bfV_{i-1}, \bfe_{k}\rangle^{2}) = O(n). 
$$
For the other term, by~\eqref{est.expo.mixing} we have
\beqnas
|\bbE(\sum^{n}_{i \neq j}m_{i, n}\langle\bfV_{i-1}, \bfe_{k}\rangle m_{j, n}\langle\bfV_{j-1}, \bfe_{k}\rangle| &=& \sum^{n}_{i \neq j}|\Cov(m_{i, n}\langle\bfV_{i-1}, \bfe_{k}, \rangle m_{j, n}\langle\bfV_{j-1}, \bfe_{k}\rangle)| \\
&\leq& \sum^{n}_{i \neq j}\omega^{|i-j|}\|m_{i, n}\|_{2}\|m_{j, n}\|_{2}, 
\eeqnas
Thus when $d =2$, we have
\beqnas
|\bbE(\sum^{n}_{i \neq j}m_{i, n}\langle\bfV_{i-1}, \bfe_{k}\rangle m_{j, n}\langle\bfV_{j-1}, \bfe_{k}\rangle| = O(\log \log n), 
\eeqnas
and when $d=3$, we get
\beqnas
|\bbE(\sum^{n}_{i \neq j}m_{i, n}\langle\bfV_{i-1}, \bfe_{k}\rangle m_{j, n}\langle\bfV_{j-1}, \bfe_{k}\rangle| = O(1).
\eeqnas
Thus we finish the proof by putting all the terms together. 
\epf

\blem
\beqna\label{m.est.tricky}
\|\sum^{n}_{i = 1}\frac{1}{d\sigma^{2}_{d}n \log n}\bbE(\tilde{\xi}^{2}_{i, n}|\bdeta) - 1\|_{2} \leq O(\frac{1}{(\log n)^{1 - \frac{\gamma}{2}}}).
\eeqna
\elem

\bpf
By definition we have
\beqnas
& &\|\sum^{n}_{i = 1}\frac{1}{d\sigma^{2}_{d}n \log n}\bbE(\tilde{\xi}^{2}_{i, n}|\bdeta) - 1\|_{2}\\
 &\leq& \frac{1}{d\sigma^{2}_{d}n \log n}\sqrt{\sum^{n}_{i = 1}\bbE (\bbE(\tilde{\xi}^{2}_{i, n}|\bdeta) - d\sigma^{2}_{d}\log n )^{2}} \\
&  & \ \ \ \ \ +  \frac{1}{d\sigma^{2}_{d}n \log n}\sqrt{\sum^{n}_{i \neq j}\bbE(\bbE(\tilde{\xi}^{2}_{i, n}|\bdeta) - d\sigma^{2}_{d}\log n)(\bbE(\tilde{\xi}^{2}_{j, n}|\bdeta) - d\sigma^{2}_{d}\log n)}.
\eeqnas
For the first term, direct computations and \eqref{m.est.tilde.2.2} yield
 \beqnas
& &\frac{1}{d\sigma^{2}_{d}n \log n}\sqrt{\bbE \sum^{n}_{i = 1}\big(\bbE((\tilde{\xi}^{2}_{i, n}|\bdeta) - d\sigma^{2}_{d}\log n \big)^{2}} \\
&\leq& \frac{1}{d\sigma^{2}_{d}n \log n}\sqrt{  \sum^{n}_{i = 1}\bbE(\bbE(\tilde{\xi}^{2}_{i, n}|\bdeta)^{2})} 
 = O(\frac{1}{(\log n)^{1 - \frac{\gamma}{2}}}).
\eeqnas

For the second term, 
by the exponential mixing estimate~\eqref{est.expo.mixing} and \eqref{m.est.tilde.2.2}, we derive that
\beqnas
& &\sum^{n}_{i \neq j}\bbE(\bbE(\tilde{\xi}^{2}_{i, n}|\bdeta) - d\sigma^{2}_{d}\log n)(\bbE(\tilde{\xi}^{2}_{j, n}|\bdeta) - d\sigma^{2}_{d}\log n)\\
&=&\sum^{n}_{i \neq j}\Cov(\bbE(\tilde{\xi}^{2}_{i, n}|\bdeta) - d\sigma^{2}_{d}\log n, \bbE(\tilde{\xi}^{2}_{j, n}|\bdeta) - d\sigma^{2}_{d}\log n) + (\bbE\tilde{\xi}^{2}_{i, n} - d\sigma^{2}_{d}\log n)(\bbE\tilde{\xi}^{2}_{j, n} - d\sigma^{2}_{d}\log n)\\
&\leq&  C\|\bbE(\tilde{\xi}^{2}_{n, n}|\bdeta) - d\sigma^{2}_{d}\log n\|^{2}_{2} +  (\bbE\tilde{\xi}^{2}_{i, n} - d\sigma^{2}_{d}\log n)(\bbE\tilde{\xi}^{2}_{j, n} - d\sigma^{2}_{d}\log n)\\
&=&  O(n(\log n)^{\gamma})+ n(n-1)O((\log \log n)^{2}), 
\eeqnas
leading to
\beqnas
\frac{1}{d\sigma^{2}_{d}n \log n}\sqrt{\sum^{n}_{i \neq j}\bbE(\bbE(\tilde{\xi}^{2}_{i, n}|\bdeta) - d\sigma^{2}_{d}\log n)(\bbE(\tilde{\xi}^{2}_{i, n}|\bdeta) - d\sigma^{2}_{d}\log n)} 
 = O(\frac{\log \log n}{ \log n}).
 \eeqnas
Finally we obtain
 \beqnas
 \|\sum^{n}_{i = 1}\frac{1}{d\sigma^{2}_{d}n \log n}\bbE(\tilde{\xi}^{2}_{i}|\bdeta) - 1\|_{2} =  O(\frac{1}{(\log n)^{1 - \frac{\gamma}{2}}})+  O(\frac{\log \log n}{ \log n}) = O(\frac{1}{(\log n)^{1 - \frac{\gamma}{2}}}).
 \eeqnas

\epf

 \section{Proof of Theorem~\ref{thm1}}
 In this section we prove our main Theorem~\ref{thm1}. The aim is to give the convergence rate of
 $$
 \sup_{h \in  \Lip(1)}|\bbE(h(\bfW_{n})) - \bbE(h(\bfZ))|.
 $$
 
Notice that  for any $h \in \Lip(1)$, we have
\beqnas
|\bbE(h(\bfW_{n})) - \bbE(h(\bfZ))| 
&\leq& \bbE |(h(\bfW_{n}) - h(\bfW'_{n})| +  |\bbE h(\bfW'_{n})  - \bbE(h(\bfZ))|.
\eeqnas

By the estimate~\eqref{estimate.wn}, we know that
\beqna\label{est1}
\nonumber && \bbE |h(\bfW_{n}) - h(\bfW'_{n})| \leq  \bbE \|D h\|\|\bfW_{n} - \bfW'_{n}\| \\
 &\leq&   \bbE \|\bfW_{n} - \bfW'_{n}\| =  O(\frac{1}{(\log n)^{\frac{\gamma + 1}{2}}}).
\eeqna

To estimate the term $ |\bbE h(\bfW'_{n})  - \bbE(h(\bfZ)) |$, we separate it into two parts,
\beqnas
|\bbE h(\bfW'_{n})  - \bbE h(\bfZ)| \leq  |\bbE h(\bfW'_{n})  -  \bbE h(\tilde{\bfW}_{n})| + |\bbE h(\tilde{\bfW}_{n})  - \bbE h(\bfZ)|,
\eeqnas
where
$$
 \tilde{\bfW}_{n} = \frac{1}{\sigma_{d}\sqrt{n \log n}}\sum^{n}_{i=1} \tilde{\xi}_{i, n}\bfV_{i-1}.
$$
Notice that
\beqnas\label{extra.term2}
& &|\bbE h(\bfW'_{n})  -  \bbE h(\tilde{\bfW}_{n})| \leq  \|\nabla h\|_{0}\bbE \| \bfW'_{n} - \tilde{\bfW}_{n}\| \\
&\leq&   \frac{1}{\sigma_{d}\sqrt{n \log n}}\sqrt{\bbE \|\sum^{n}_{i=1}m_{i, n}\bfV_{i-1}\|^{2}},
\eeqnas
thus by \eqref{est.sum.min.d2} and \eqref{est.sum.min.d3}, we have 
\beqna\label{est2.d2}
|\bbE h(\bfW'_{n})  -  \bbE h(\tilde{\bfW}_{n})| = O(\sqrt{\frac{\log \log n}{\log n}}),
\eeqna
when $d=2$, and
\beqna\label{est2.d3}
|\bbE h(\bfW'_{n})  -  \bbE h(\tilde{\bfW}_{n})| = O(\frac{1}{\sqrt{\log n}}).
\eeqna
when $d \geq 3$.

As for the term $ |\bbE h(\tilde{\bfW}_{n})  - \bbE h(\bfZ)|$, we apply Stein's method to obtain its convergence rate. 

 \subsection{The estimate of $|\bbE h(\tilde{\bfW}_{n})  - \bbE h(\bfZ)|$}
 
In this section we use the approach of "exchangeable pair" in Stein's method to estimate  $|\bbE h(\tilde{\bfW}_{n})  - \bbE h(\bfZ)|$. The main obstacle here is that $\{(\xi_{n}, \bdeta_{n})\}$ are not independent, and the estimates of their moments are quite delicate. We follow the idea from \cite{M-T16}, i.e. to make use of the fact that $\{\xi_{i}\}$ are independent  under the conditional expectation $\bbE(\cdot |\bdeta)$, although no longer identically distributed. 

Thus to apply Stein's method, we need to adapt the technique under the conditional expectation $\bbE(\cdot |\bdeta)$. Inspired by Chatterjee-Meckes \cite{C-M07}, we first construct the exchangeable pair $( \tilde{\bfW}_{n}, \tilde{\bfW}'_{n})$ under $\bbE(\cdot |\bdeta)$. 
%


For each $i$,  let $\tilde{\xi}'_{i}$ be the random variable on the same probability space as  $\tilde{\xi}_{i}$, which is identically distributed as $\tilde{\xi}_{i}$, but independent of $\{\tilde{\xi}_{i}\}$ under $\bbE(\cdot | \bdeta)$. 
Moreover, let $I$ an index number random  variable  taking values  on $\{1, 2, \dots, n\}$  with $\bbP(I = i) = \frac{1}{n}$.

Then we define 
\beqnas
\tilde{\bfW}'_{n} &:=& \tilde{\bfW}_{n} - \frac{1}{\sigma_{d}\sqrt{n \log n}} \tilde{\xi}_{I,n}\bfV_{I-1} +  \frac{1}{\sigma_{d}\sqrt{n \log n}} \tilde{\xi}'_{I, n}\bfV_{I-1}\\
&=& \tilde{\bfW}_{n} - \frac{1}{\sigma_{d}\sqrt{n \log n}} \bfV_{I-1}(\tilde{\xi}_{I, n} - \tilde{\xi}'_{I, n}), 
\eeqnas
i.e.
\beqnas
\tilde{\bfW}'_{n} -  \tilde{\bfW}_{n}  =  \frac{1}{\sigma_{d}\sqrt{n \log n}} \bfV_{I-1}(\tilde{\xi}'_{I, n} - \tilde{\xi}_{I, n}).
\eeqnas

Recall the definition~\ref{def.exc.pair} and the conditions~\eqref{exc.pair.eq1} and \eqref{exc.pair.eq2}, we verify the following proposition, where \eqref{exchange.c2} and \eqref{exchange.c3}  are analogous to \eqref{exc.pair.eq1} and \eqref{exc.pair.eq2}.
\bprop
$(\tilde{\bfW}_{n}, \tilde{\bfW}'_{n})$ is an exchangeable pair under $\bbE(\cdot |\bdeta)$.  Moreover, we have
\beqna\label{exchange.c2}
\bbE(\bbE(\tilde{\bfW}'_{n} - \tilde{\bfW}_{n} |\tilde{\bfW}_{n})|\bdeta) =  - \frac{1}{n}\bbE(\tilde{\bfW}_{n}|\bdeta), 
\eeqna
and 
\beqna\label{exchange.c3}
& &\bbE( (\tilde{\bfW}'_{n} -  \tilde{\bfW}_{n})(\tilde{\bfW}'_{n} -  \tilde{\bfW}_{n})^{t}|\sigma(\tilde{\bfW}_{n}))\\
\nonumber &=&  \frac{1}{\sigma^{2}_{d}n^{2} \log n}\sum^{n}_{i=1}\bfV_{i-1}\bfV^{t}_{i-1}( \bbE( (\tilde{\xi}_{i, n})^{2}|\bdeta)  + \tilde{\xi}^{2}_{i, n}).
\eeqna
\eprop

\bpf
First we prove that under $\bbE(\cdot | \bdeta)$, $(\tilde{\bfW}_{n}, \tilde{\bfW}'_{n})$ has the same distribution as $(\tilde{\bfW}'_{n}, \tilde{\bfW}_{n})$, i.e. for any bounded continuous function $f$, 
\beqnas
\bbE ( f(\tilde{\bfW}_{n}, \tilde{\bfW}'_{n})| \bdeta) = \bbE(f(\tilde{\bfW}'_{n}, \tilde{\bfW}_{n}) | \bdeta),
\eeqnas
which leads to
\beqna\label{exchange.c1}
\bbE f(\tilde{\bfW}_{n}, \tilde{\bfW}'_{n}) = \bbE f(\tilde{\bfW}'_{n}, \tilde{\bfW}_{n}).
\eeqna
To see this, notice that
\beqnas
& &\bbE ( f(\tilde{\bfW}_{n}, \tilde{\bfW}'_{n})| \bdeta) \\
&=& \frac{1}{n}\sum^{n}_{i= 1}\bbE(\bbE ( f(\tilde{\bfW}'_{n}+  \frac{1}{\sigma_{d}\sqrt{n \log n}} \bfV_{i-1}(\tilde{\xi}_{i, n} - \tilde{\xi}'_{i, n}), \tilde{\bfW}_{n} - \frac{1}{\sigma_{d}\sqrt{n \log n}} \bfV_{i-1}(\tilde{\xi}_{i, n} - \tilde{\xi}'_{i, n}))| I = i)| \bdeta)\\
&=& \bbE(f(\tilde{\bfW}'_{n}, \tilde{\bfW}_{n})| \bdeta), 
\eeqnas 
where the last line is due to the fact that under $\bbE(\cdot |\bdeta)$, $\{\tilde{\xi}_{i, n}, \tilde{\xi}'_{i, n}\}$ are independent and for each $i$, $\tilde{\xi}_{i, n}, \tilde{\xi}'_{i, n}$ are identically distributed. 

Moreover, we have
\beqnas
& &\bbE(\bbE(\tilde{\bfW}'_{n}|\tilde{\bfW}_{n})|\bdeta) \\
& =& \sum^{n}_{i= 1}\bbE(\bbE ( \tilde{\bfW}_{n} - \frac{1}{\sigma_{d}\sqrt{n \log n}} \bfV_{i-1}(\tilde{\xi}_{i, n} - \tilde{\xi}'_{i, n}) |\tilde{\bfW}_{n},  I = i)|\bdeta) \\
& =&  \bbE(\tilde{\bfW}_{n}|\bdeta)  - \frac{1}{\sigma_{d}\sqrt{n \log n}}  \sum^{n}_{i= 1}\bbE(\bbE(\bfV_{i-1}(\tilde{\xi}_{i, n} - \tilde{\xi}'_{i, n})|\tilde{\bfW}_{n},  I = i)|\bdeta).
\eeqnas
%
The independence of $\tilde{\xi}'_{i, n}$ and $\tilde{\xi}_{i, n}$ under $\bbE(\cdot |\bdeta)$ implies that
\beqnas
\bbE(\bbE(\bfV_{i-1}\tilde{\xi}'_{i, n}|\tilde{\bfW}_{n})|\bdeta)  = 0,
\eeqnas
leading to
\beqnas
\bbE(\bbE(\tilde{\bfW}'_{n} - \tilde{\bfW}_{n} |\tilde{\bfW}_{n})|\bdeta) &=&  - \frac{1}{\sigma_{d}\sqrt{n \log n}}  \sum^{n}_{i= 1}\bbE(\bbE(\bfV_{i-1}\tilde{\xi}_{i, n}|\tilde{\bfW}_{n},  I = i)|\bdeta)\\
&=&  - \frac{1}{\sigma_{d}\sqrt{n \log n}}  \sum^{n}_{i= 1}\bbE(\bbE(\bfV_{i-1}\tilde{\xi}_{i, n}|\tilde{\bfW}_{n}|\bdeta)\bbP(I = i)\\
&=& - \frac{1}{n}\bbE(\tilde{\bfW}_{n}|\bdeta),
\eeqnas
in that  $I$ is uniformly distributed and independent of $\{\tilde{\xi}_{i, n}, \tilde{\xi}'_{i, n}\}$.

It can be verified in the same way that 
\beqnas
& &\bbE( (\tilde{\bfW}'_{n} -  \tilde{\bfW}_{n})(\tilde{\bfW}'_{n} -  \tilde{\bfW}_{n})^{t}|\sigma(\tilde{\bfW}_{n}))\\
\nonumber&=&  \frac{1}{\sigma^{2}_{d}n \log n} \sum^{n}_{i=1}\bbE( (\tilde{\xi}'_{i, n} - \tilde{\xi}_{i, n})^{2}\bfV_{i-1}\bfV^{t}_{i-1}|\sigma(\xi_{j}, \bdeta_{j}),  I = i)\bbP(I = i)\\
\nonumber &=&  \frac{1}{\sigma^{2}_{d}n^{2} \log n}\sum^{n}_{i=1}\bfV_{i-1}\bfV^{t}_{i-1}( \bbE( (\tilde{\xi}_{i, n})^{2}|\bdeta)  + \tilde{\xi}^{2}_{i, n}).
\eeqnas

\epf
%

Now we turn to the estimate of $|\bbE h(\tilde{\bfW}_{n})  - \bbE h(\bfZ)|$, using the exchangeable pair  $(\tilde{\bfW}_{n}, \tilde{\bfW}'_{n})$.

Given $h \in \Lip(1)$, let $f_{h}$ be the solution to the Stein's equation in multi-dimension~\eqref{stein.eq.multid}. Following the discussion in Section 2.2 and taking 
$$
f(\bfx, \bfx') = \langle (\bfx' -  \bfx), \nabla f_{h}(\bfx') +  \nabla f_{h}(\bfx) \rangle
$$ in \eqref{exchange.c1},  we derive that
\beqna\label{eq.R}
0 &=& \frac{n}{2}\bbE[\langle (\tilde{\bfW}'_{n} -  \tilde{\bfW}_{n}), \nabla f_{h}(\tilde{\bfW}'_{n}) +  \nabla f_{h}( \tilde{\bfW}_{n}) \rangle]\\
\nonumber &=& \frac{n}{2}\bbE[\langle (\tilde{\bfW}'_{n} -  \tilde{\bfW}_{n}), \nabla f_{h}(\tilde{\bfW}'_{n}) - \nabla f_{h}( \tilde{\bfW}_{n}) \rangle + 2\langle (\tilde{\bfW}'_{n} -  \tilde{\bfW}_{n}), \nabla f_{h}( \tilde{\bfW}_{n}) \rangle]\\
\nonumber &=& \frac{n}{2}\bbE[\langle (\tilde{\bfW}'_{n} -  \tilde{\bfW}_{n})(\tilde{\bfW}'_{n} -  \tilde{\bfW}_{n})^{t}, {\rm Hess} f_{h}(\tilde{\bfW}_{n})\rangle + 2\langle (\tilde{\bfW}'_{n} -  \tilde{\bfW}_{n}), \nabla f_{h}( \tilde{\bfW}_{n}) \rangle + R],
\eeqna
where $R$ is given by
\beqnas
R  &=&  \frac{1}{2} \sum^{d}_{j, k, l=1}D^{3}_{jkl}f_{h}(\tilde{\bfW}_{n} + \tau(\tilde{\bfW}'_{n} -  \tilde{\bfW}_{n}))(\tilde{\bfW}'_{n} -  \tilde{\bfW}_{n})_{j}(\tilde{\bfW}'_{n} -  \tilde{\bfW}_{n})_{k}(\tilde{\bfW}'_{n} -  \tilde{\bfW}_{n})_{l}.
\eeqnas


Then by  \eqref{exchange.c2} and \eqref{exchange.c3} we have
\beqnas
0 &=& \bbE[\langle \frac{1}{2\sigma^{2}_{d}n \log n}\sum^{n}_{i=1}\bfV_{i-1}\bfV^{t}_{i-1}( \bbE(\tilde{\xi}_{i, n}^{2}|\bdeta)  + \tilde{\xi}^{2}_{i, n}), {\rm Hess} f_{h}(\tilde{\bfW}_{n})\rangle - \langle \tilde{\bfW}_{n}, \nabla f_{h}( \tilde{\bfW}_{n}) \rangle + \frac{n}{2}R],
\eeqnas
such that  the Stein's equation~\eqref{stein.eq.multid} leads to 
\beqnas
& &\bbE h(\tilde{\bfW}_{n}) - \bbE h(\bfZ) \\
&=& \bbE [\tr D^{2}f_{h}( \tilde{\bfW}_{n}) - \langle \tilde{\bfW}_{n}, \nabla f_{h}( \tilde{\bfW}_{n}) \rangle]\\
 &=& \bbE[\langle \big(\Id - \frac{1}{2\sigma^{2}_{d}n \log n}\sum^{n}_{i=1}\bfV_{i-1}\bfV^{t}_{i-1}( \bbE(\tilde{\xi}_{i, n}^{2}|\bdeta)  + \tilde{\xi}^{2}_{i, n})\big), {\rm Hess} f_{h}(\tilde{\bfW}_{n})\rangle  - \frac{n}{2}R].
\eeqnas

Notice that by defining $\bbE (D^{3}_{jkl}f_{h}(\tilde{\bfW}_{n} + \tau(\tilde{\bfW}'_{n} -  \tilde{\bfW}_{n}))(\tilde{\xi}'_{i, n} - \tilde{\xi}_{i, n})^{3}|\bdalpha_{i}) = \tilde{f}_{3}(\bdalpha_{i})$, we have
\beqnas
\bbE R &=&\frac{1}{2n\sigma^{3}_{d}(n \log n)^{\frac{3}{2}}}\\
& & \hspace{1cm} \cdot \sum^{n}_{i=1} \sum^{d}_{j, k, l=1}\bbE(V_{i-1, j}V_{i-1, k}V_{i-1, l}\bbE (D^{3}_{jkl}f_{h}(\tilde{\bfW}_{n} + \tau(\tilde{\bfW}'_{n} -  \tilde{\bfW}_{n}))(\tilde{\xi}'_{i, n} - \tilde{\xi}_{i, n})^{3}|\bdeta))\\
 &=&\frac{1}{2n\sigma^{3}_{d}(n \log n)^{\frac{3}{2}}}\sum^{n}_{i=1} \sum^{d}_{j, k, l=1}\bbE(V_{i-1, j}V_{i-1, k}V_{i-1, l}\tilde{f}_{3}(\bdalpha_{i})) = 0,
\eeqnas
where the last line is due to~\eqref{estimate.trois}. 

Now we give the estimates of the first term. We first assume that $h \in \Lip(1)$ with  up to the second bounded derivatives. 
\blem\label{lm4.2}
\beqna\label{w1.est.1}
\nonumber & &| \bbE \langle \big(\Id - \frac{1}{2\sigma^{2}_{d}n \log n}\sum^{n}_{i=1}\bfV_{i-1}\bfV^{t}_{i-1}( \bbE(\tilde{\xi}_{i, n}^{2}|\bdeta)  + \tilde{\xi}^{2}_{i, n})\big), {\rm Hess} f_{h}(\tilde{\bfW}_{n})\rangle| \\
&=& \sup \|D^{2}f_{h}\|O(\frac{1}{(\log n)^{1 - \frac{\gamma}{2}}}).
\eeqna
\elem

\bpf
We first write 
\beqnas
& &\bbE \langle \big(\Id - \frac{1}{2\sigma^{2}_{d}n \log n}\sum^{n}_{i=1}\bfV_{i-1}\bfV^{t}_{i-1}(  \bbE(\tilde{\xi}_{i, n}^{2}|\bdeta)  + \tilde{\xi}^{2}_{i, n})\big), {\rm Hess} f_{h}(\tilde{\bfW}_{n})\rangle\\
&=& \sum^{d}_{j =1}\bbE [\big(1 - \frac{1}{2\sigma^{2}_{d}n \log n}\sum^{n}_{i=1}V^{2}_{i-1, j}(\bbE(\tilde{\xi}_{i, n}^{2}|\bdeta)  + \tilde{\xi}^{2}_{i, n})\big)D^{2}_{jj}f_{h}(\tilde{\bfW}_{n})]\\
& & -  \frac{1}{2\sigma^{2}_{d}n \log n} \sum^{d}_{j \neq k} \bbE \big(\sum^{n}_{i=1}V_{i-1, j}V_{i-1, k}(\bbE(\tilde{\xi}_{i, n}^{2}|\bdeta)  + \tilde{\xi}^{2}_{i, n})D^{2}_{jk}f_{h}(\tilde{\bfW}_{n})\big)\\
&:=& I_{1} + I_{2}.
\eeqnas
%
As for $I_{1}$,   we have
\beqnas
|I_{1}| & \leq& \sum^{d}_{j =1}\bbE [\big(1 - \frac{1}{2d\sigma^{2}_{d}n \log n}\sum^{n}_{i=1}(\bbE(\tilde{\xi}_{i, n}^{2}|\bdeta)  + \tilde{\xi}^{2}_{i, n})\big)D^{2}_{jj}f_{h}(\tilde{\bfW}_{n})]\\
 & & \ \ \ \ +  \frac{1}{2\sigma^{2}_{d}n \log n}\sum^{d}_{j =1}|\bbE \big(\sum^{n}_{i = 1}(V^{2}_{i-1, j} - \frac{1}{d})(\bbE( (\bbE( \tilde{\xi}_{i, n}^{2}|\bdeta)  + \tilde{\xi}^{2}_{i, n})D^{2}_{jj}f_{h}(\tilde{\bfW}_{n})|\bdeta)\big)|\\
 &:=& I_{1, 1}  + I_{1, 2}.
  \eeqnas

For the first term $I_{1, 1}$, by \eqref{m.est.tricky} we derive 
\beqnas
|I_{1, 1}| &\leq& \sum^{d}_{j =1}|\bbE [\big(1 - \frac{1}{2d\sigma^{2}_{d}n \log n}\sum^{n}_{i=1}(\bbE(\tilde{\xi}_{i, n}^{2}|\bdeta)  + \tilde{\xi}^{2}_{i, n})\big)D^{2}_{jj}f_{h}(\tilde{\bfW}_{n})]|\\
&\leq& \sup\|D^{2}f_{h}\| \| 1- \sum^{n}_{i = 1}\frac{1}{d\sigma^{2}_{d}n \log n}\bbE(\tilde{\xi}^{2}_{i, n}|\bdeta)\|_{2} =   \sup\|D^{2}f_{h}\| O(\frac{1}{(\log n)^{1 - \frac{\gamma}{2}}}).
\eeqnas

The second term $I_{1, 2}$ is estimated by~\eqref{estimate.deux2}.
\beqnas
|I_{1, 2}| &\leq& \frac{1}{2\sigma^{2}_{d}n \log n}|\bbE \big(\sum^{n}_{i = 1}(\bfV^{2}_{i-1, j} - \frac{1}{d})\bbE( \bbE( \tilde{\xi}_{i, n}^{2}|\bdeta)  + \tilde{\xi}^{2}_{i, n})D^{2}_{jj}f_{h}(\tilde{\bfW}_{n})|\bdeta)| \\
&\leq&  \frac{C_{d} \sup\|D^{2}f_{h}\| }{n \log n}\|\bbE( \tilde{\xi}_{i, n}^{2}|\bdeta)\|_{2} =    \sup\|D^{2}f_{h}\| O(\frac{1}{n^{\frac{1}{2}}(\log n)^{1 -\frac{\gamma}{2} }}).
\eeqnas

Therefore we obtain the estimate of $I_{1}$, 
\beqna\label{l42.est.1}
|I_{1}| \leq \sup\|D^{2}f_{h}\| O(\frac{1}{(\log n)^{1 - \frac{\gamma}{2}}}).
\eeqna

By \eqref{estimate.deux1} and \eqref{m.est.tilde.2.2} we derive that
\beqna\label{l42.est.2}
\nonumber |I_{2}|  &\leq& \frac{1}{2\sigma^{2}_{d}n \log n} \sum^{d}_{j \neq k} \sum^{n}_{i=1}|\bbE \big(V_{i-1, j}V_{i-1, k}\bbE((\bbE(\tilde{\xi}_{i, n}^{2}|\bdeta)  + \tilde{\xi}^{2}_{i, n})D^{2}_{jk}f_{h}(\tilde{\bfW}_{n})|\bdeta)\big)|\\
\nonumber &\leq&\frac{C_{d}}{n \log n} \sum^{n}_{i=1}|\bbE \big(V_{i-1, j}V_{i-1, k}( \bbE(\tilde{\xi}_{i, n}^{2}|\bdeta)\bbE(D^{2}_{jk}f_{h}(\tilde{\bfW}_{n})|\bdeta) + \bbE(\tilde{\xi}_{i, n}^{2}D^{2}_{jk}f_{h}(\tilde{\bfW}_{n})|\bdeta))\big)|\\
 &\leq& \frac{C_{d}}{n \log n} \sup\|D^{2}f_{h}\| \|\bbE(\tilde{\xi}_{i, n}^{2}|\bdeta)\|_{2} =    \sup\|D^{2}f_{h}\| O(\frac{1}{n^{\frac{1}{2}} (\log n)^{1 - \frac{\gamma}{2}}}).
\eeqna

Combining \eqref{l42.est.1} and  \eqref{l42.est.2} together leads to the estimate~\eqref{w1.est.1}.
\epf

Thus with~\eqref{w1.est.1} we yield 
\beqna\label{est3}
|\bbE h(\tilde{\bfW}_{n}) - \bbE h(\bfZ)| \leq  \sup \|D^{2}f_{h}\|O(\frac{1}{(\log n)^{1 - \frac{\gamma}{2}}}).
\eeqna

%
%

Now we are in the position to prove our main Theorem~\ref{thm1}. 

\bpf
As in Chatterjee-Meckes~\cite{C-M07}, for $h \in \Lip(1)$ we first consider the smooth function $h_{\epsilon} = h \ast \phi_{\epsilon}$, where $\phi_{\epsilon}$ is the density of $\cN(0, \epsilon^{2}\Id_{d})$, such that we have
\benum
{\item
\beqna\label{bound.sm}
\|h_{\epsilon} - h\|_{\infty} \rightarrow 0,
\eeqna
}
{\item
\beqna\label{bound.smderi.k}
\|D^{k}(h_{\epsilon})\| \leq \|D^{k}(h)\|.
\eeqna
}
\eenum

With \eqref{est1}, \eqref{est2.d2} (or  \eqref{est2.d3}), \eqref{est3}, we derive that when $d=2$
\beqna\label{estimate.mt2}
& &|\bbE h_{\epsilon}(\bfW_{n}) - \bbE h_{\epsilon}(\bfZ) |\\
\nonumber &\leq& O(\frac{1}{(\log n)^{\frac{\gamma + 1}{2}}}) + O(\sqrt{\frac{\log \log n}{\log n}}) + \sup \|D^{2}f_{h_{\epsilon}}\|O(\frac{1}{(\log n)^{1 - \frac{\gamma}{2}}}), 
\eeqna
while for $d \geq 3$, 
\beqna\label{estimate.mt3}
& &|\bbE h_{\epsilon}(\bfW_{n}) - \bbE h_{\epsilon}(\bfZ) |\\
\nonumber &\leq& O(\frac{1}{(\log n)^{\frac{\gamma + 1}{2}}}) + O(\frac{1}{\sqrt{\log n}}) +\sup \|D^{2}f_{h_{\epsilon}}\|O(\frac{1}{(\log n)^{1 - \frac{\gamma}{2}}}) 
\eeqna

Therefore, by letting $\epsilon \rightarrow 0$ and with~\eqref{bound.deri2}, \eqref{bound.smderi.k},  we obtain the estimates~\eqref{thm1.est1} and ~\eqref{thm1.est2}. 


%
%
\epf

As a corollary of Theorem~\ref{thm1}, we immediately get the Berry-Essen type estimate. See Erickson~\cite{Erickson74}, Chen-Shao~\cite{C-M07} for the connection between the Berry-Essen type bounds and the bound with respect to Lipschitz functions. We point out that Proposition~\ref{thm2} can also be obtained by the Stein's method, but we fail to improve the estimates. The computations are tedious so we omit them here.

\bprop\label{thm2}
Let $\bfz \in \bbR^{d}$. Under assumptions $\bfA$ on the initial data, we have when $d=2$,
\beqna
|\bbP(\bfW_{n} \leq \bfz) - \Phi(\bfz)| =   O(\frac{\log \log n}{\log n})^{\frac{1}{4}}),
\eeqna
and when $d \geq 3$
\beqna
|\bbP(\bfW_{n} \leq \bfz) -  \Phi(\bfz)| =   O(\frac{1}{(\log n)^{\frac{1}{4}}}),
\eeqna
where $\Phi$ is the Gaussian distribution function on $\bbR^{d}$.
\eprop

 \section{Proof of Theorem~\ref{thm3}}
 In this section we prove our main  Theorem~\ref{thm3}, which gives the convergence rate of the continuous time displacement $\bfX_{t}$.  
Let 
\beqnas
n_{t} = [\bar{\xi}^{-1}t],
\eeqnas
where $\bar{\xi} = \frac{1}{v_{d-1}}$ is the mean free path length. Also define
$$
\bfW_{t} =  \frac{1}{\Sigma_{d}\sqrt{t\log t}}\bfX_{t}.
$$
To estimate the convergence rate of the continuous process $\bfW_{t}$, we compare it with the discrete case $\bfW_{n_{t}}$. For $\bfz \in \bbR^{d}$, let $\bfz_{1} = z_{1}\frac{\bfz}{\|\bfz\|}$ with $z_{1}$ to be decided later. Notice that
\beqna\label{ineq1}
\nonumber & &|\bbP(\bfW_{t} \geq \bfz)  -  \bbP(\bfZ \geq \bfz) | \\
&\leq& |\bbP(\bfW_{t} - \bfW_{n_{t}} \geq \bfz_{1}) | + |\bbP(\bfW_{n_{t}} \geq \bfz - \bfz_{1})  -  \bbP(\bfZ \geq \bfz) |,
\eeqna
where  $ |\bbP(\bfW_{n_{t}} \geq \bfz - \bfz_{1})  -  \bbP(\bfZ \geq \bfz) |$ is estimated in Corollary~\ref{thm2}, i.e. for $d \geq 3$, we have
\beqna\label{ineq2.1}
\nonumber & &|\bbP(\bfW_{n_{t}} \geq \bfz - \bfz_{1})  - \bbP(\bfZ \geq \bfz) | \\
\nonumber&\leq & |\bbP(\bfW_{n_{t}} \geq \bfz - \bfz_{1})  -  \bbP(\bfZ \geq \bfz- \bfz_{1})  | + | \bbP(\bfZ \geq \bfz- \bfz_{1}) -  \bbP(\bfZ \geq \bfz) |\\
 &\leq& O(\frac{1}{(\log t)^{\frac{1}{4}}})+ O(z_{1}),
\eeqna
and for $d =2$ we have
\beqna\label{ineq2.2}
|\bbP(\bfW_{n_{t}} \geq \bfz - \bfz_{1})  - \bbP(\bfZ \geq \bfz) | \leq O(\frac{(\log \log t)^{\frac{1}{4}}}{(\log t)^{\frac{1}{4}}})+ O(z_{1}).
\eeqna

 In the following we will focus on the estimate of $|\bbP(\bfW_{t} - \bfW_{n_{t}} \geq \bfz_{1})|$. 

By definition 
\beqnas
& &\bfW_{t} -\bfW_{n_{t}} =  \frac{1}{\Sigma_{d}\sqrt{t\log t}}\big(\bfX_{t} - \frac{\sqrt{t\log t}}{\bar{\xi}^{\frac{1}{2}}\sqrt{n_{t}\log n_{t}}} \bfQ_{n_{t}}\big)\\
&=&  \frac{1}{\Sigma_{d}\sqrt{t\log t}}\big( (t - \tau_{\nu_{t}})\bfV_{\nu_{t}} + \bfQ_{\nu_{t}} - \bfQ_{n_{t}}  +(1- \frac{\sqrt{t\log t}}{\bar{\xi}^{\frac{1}{2}}\sqrt{n_{t}\log n_{t}}} )\bfQ_{n_{t}} \big), 
\eeqnas
we have
\beqna\label{ineq3}
& &|\bbP(\bfW_{t} - \bfW_{n_{t}} \geq \bfz_{1})| \\
\nonumber &\leq& \bbP( \frac{1}{\Sigma_{d}\sqrt{t\log t}}\|(t - \tau_{\nu_{t}})\bfV_{\nu_{t}}\| \geq \frac{z_{1}}{3} ) + \bbP( \frac{1}{\Sigma_{d}\sqrt{t\log t}}\| \bfQ_{\nu_{t}} - \bfQ_{n_{t}}\| \geq \frac{z_{1}}{3} )\\
\nonumber & & + \bbP(|\frac{1}{\Sigma_{d}\sqrt{t\log t}}- \frac{1}{\Sigma_{d}\bar{\xi}^{\frac{1}{2}}\sqrt{n_{t}\log n_{t}}}|\|\bfQ_{n_{t}}\| \geq \frac{z_{1}}{3})\\
\nonumber &:=& I_{1} + I_{2} + I_{3}. 
\eeqna

The first term $I_{1}$ is easy to estimate
\beqna\label{est.cont.i1}
I_{1} \leq \frac{C}{z_{1}\sqrt{t\log t}}\bbE|\xi_{\nu_{t}+1}| = O(\frac{1}{z_{1}\sqrt{t\log t}}). 
\eeqna

In order to estimate $I_{2}$ and $I_{3}$, we prove the following lemmas on $\bbE\|\bfQ_{n}\|$ and $\bbE\sup_{t \in [0, 1]}\|\bfQ_{[nt]}\|$. 

\blem
\beqna\label{est.qn}
\bbE \|\bfQ_{n}\| = O(\sqrt{n\log n}).
\eeqna
\elem
\bpf
To estimate $\bfQ_{n}$, we first make use of the truncation. Recall~\eqref{estimate.qn}, we have
\beqnas
\bbE \|\bfQ_{n}\|  \leq \bbE \|\bfQ_{n} - \bfQ'_{n}\| + \bbE \|\bfQ'_{n}\|  =   O( \frac{n^{\frac{1}{2}}}{(\log n)^{\frac{\gamma}{2}}}) + \bbE \|\bfQ'_{n}\|.
\eeqnas
As for $\bbE \|\bfQ'_{n}\|$, we have
\beqnas
& & \bbE \|\bfQ'_{n}\|^{2} \leq \bbE \|\sum^{n}_{j=1}\xi'_{j, n}\bfV_{j-1}\|^{2} \leq  \bbE \sum^{d}_{k=1}| \sum^{n}_{j=1}  \xi'_{j, n}\langle \bfe_{k}, \bfV_{j-1}\rangle|^{2}\\
  &\leq&   C_{d}\sum^{d}_{k=1} \big(\sum^{n}_{j=1}\bbE  (\xi'_{j, n})^{2}\langle \bfe_{k}, \bfV_{j-1}\rangle^{2} + \sum^{n}_{i \neq j}\bbE  \xi'_{i, n}\xi'_{j, n}\langle \bfe_{k}, \bfV_{i-1}\rangle\langle \bfe_{k}, \bfV_{j-1}\rangle \big).
\eeqnas
By~\eqref{estimate.deux2}  we derive that
\beqnas
\sum^{d}_{k=1}\sum^{n}_{j=1}|\bbE(\xi'^{2}_{j, n}\langle \bfe_{k}, \bfV_{j-1}\rangle^{2})| &=& \sum^{d}_{k=1}\sum^{n}_{j=1}|\bbE(\langle \bfe_{k}, \bfV_{j-1}\rangle ^{2} - \frac{1}{d})\bbE(\xi'^{2}_{j, n} | \bdeta))| +  \sum^{n}_{j=1}\bbE \xi'^{2}_{j, n}\\
& =&   O(n\log n) + O(r^{2}_{n, n}) = O(n\log n).
\eeqnas
Also \eqref{est.expo.mixing} yields
\beqnas
\sum^{n}_{i \neq j}|\bbE  \xi'_{i, n}\xi'_{j, n}\langle \bfe_{k}, \bfV_{i-1}\rangle\langle \bfe_{k}, \bfV_{j-1}\rangle| 
&\leq& \sum^{n}_{i \neq j} \omega^{|i-j|}\|m_{i, n}\|_{2}\|m_{j, n}\|_{2}, 
\eeqnas
thus when $d=2$, we conclude
$$
\sum^{n}_{i \neq j}|\bbE  \xi'_{i, n}\xi'_{j, n}\langle \bfe_{k}, \bfV_{i-1}\rangle\langle \bfe_{k}, \bfV_{j-1}\rangle| = O(\log \log n), 
$$
and when $d =3$
$$
\sum^{n}_{i \neq j}|\bbE  \xi'_{i, n}\xi'_{j, n}\langle \bfe_{k}, \bfV_{i-1}\rangle\langle \bfe_{k}, \bfV_{j-1}\rangle| = O(1),
$$
leading to 
\beqnas
\bbE \|\bfQ'_{n}\|  =  O(\sqrt{n \log n}), 
\eeqnas
thus finishing the proof. 
\epf

\blem
\beqna\label{est.qn.cont}
\bbE \sup_{t \in [0, 1]} \|\bfQ_{[nt]}\| = O(\sqrt{n\log n}).
\eeqna
\elem
\bpf
As in the previous lemma, by \eqref{estimate.qn} we have
\beqnas
\bbE \sup_{t  \in [0, 1]} \|\bfQ_{[nt]}\| 
& \leq&  \sum^{n}_{i=1} \bbE |\xi_{i} - \xi'_{i, n}| + \bbE \sup_{t  \in [0, 1]} \|\sum^{[nt]}_{i=1}\xi'_{i, n}\bfV_{i-1}\|\\
&\leq& O( \frac{n^{\frac{1}{2}}}{(\log n)^{\frac{\gamma}{2}}}) + \bbE \sup_{t  \in [0, 1]} \|\sum^{[nt]}_{i=1}\xi'_{i, n}\bfV_{i-1}\|.
\eeqnas

For the second term, we separate it into two parts, 
\beqnas
 \bbE \sup_{t} \|\sum^{[nt]}_{i=1}\xi'_{i, n}\bfV_{i-1}\| \leq  \bbE \sup_{t} \|\sum^{[nt]}_{i=1}m_{i, n}\bfV_{i-1}\| + \bbE \sup_{t} \|\sum^{[nt]}_{i=1}\tilde{\xi}_{i, n}\bfV_{i-1}\|.
\eeqnas

Recall Lemma 13.4 in \cite{M-T16}, for any $\beta \in \bbR$ we have
\beqnas
\bbP(\frac{1}{\sqrt{n\log n}}\sup_{t} \|\sum^{[nt]}_{i=1}m_{i, n}\bfV_{i-1}\| \geq \beta) \leq O(\frac{1}{\beta^{2}\log n}), 
\eeqnas
such that
\beqnas
& &\frac{1}{\sqrt{n\log n}}\bbE \sup_{t} \|\sum^{[nt]}_{i=1}m_{i, n}\bfV_{i-1}\| \\
&\leq& \frac{1}{\sqrt{n\log n}}\bbE (\sup_{t} \|\sum^{[nt]}_{i=1}m_{i, n}\bfV_{i-1}\|, \frac{1}{\sqrt{n\log n}}\sup_{t} \|\sum^{[nt]}_{i=1}m_{i, n}\bfV_{i-1}\|  \leq \frac{1}{\sqrt{\log n}}) \\
&  & \ \ \ \  + \frac{1}{\sqrt{n\log n}}\bbE (\sup_{t} \|\sum^{[nt]}_{i=1}m_{i, n}\bfV_{i-1}\|, \frac{1}{\sqrt{n\log n}}\sup_{t} \|\sum^{[nt]}_{i=1}m_{i, n}\bfV_{i-1}\|  \geq \frac{1}{\sqrt{\log n}}) \\
&\leq&  \frac{1}{\sqrt{\log n}} + \int^{\infty}_{\frac{1}{\sqrt{\log n}}}\bbP(\frac{1}{\sqrt{n\log n}}\sup_{t} \|\sum^{[nt]}_{i=1}m_{i, n}\bfV_{i-1}\| \geq u)du \\
&\leq&  \frac{1}{\sqrt{\log n}} +O(\frac{1}{\sqrt{\log n}}) = O(\frac{1}{\sqrt{\log n}}).
\eeqnas

Now we turn the last term $ \bbE \sup_{t} \|\sum^{[nt]}_{i=1}\tilde{\xi}_{i, n}\bfV_{i-1}\|$. We apply the maximal inequality (see Billingsley~\cite{Billing68}), as mentioned in \cite{M-T16}, which leads to
\beqnas
& &\bbP( \frac{1}{\sqrt{n\log n}}\sup_{t \in [0, 1]}\|\sum^{[nt]}_{i=1}\tilde{\xi}_{i, n}\bfV_{i-1}\| \geq \beta) = \bbE(\bbP( \frac{1}{\sqrt{n\log n}}\sup_{t \in [0, 1]}\|\sum^{[nt]}_{i=1}\tilde{\xi}_{i,n}\bfV_{i-1}\| \geq \beta |\bdeta))\\
&\leq&  2\bbE(\bbP(\sup_{t \in [0, 1]}\|\sum^{[nt]}_{i=1}\tilde{\xi}_{i, n}\bfV_{i-1}\| \geq \beta\sqrt{n\log n} - \sqrt{2\cA^{2}_{[nt]}}|\bdeta))\\
&\leq&  2\bbP(\|\sum^{[nt]}_{i=1}\tilde{\xi}_{i, n}\bfV_{i-1}\| \geq (\beta - K \sqrt{\delta})\sqrt{n\log n}) + 2\bbP(\cA^{2}_{[nt]} > K^{2}\delta n \log n),
\eeqnas
where $\cA^{2}_{[nt]} =\bbE((\sum^{[nt]}_{i=1}\tilde{\xi}_{i,n}\bfV_{i-1})^{2}|\bdeta) =  \sum^{[nt]}_{i=1}\bbE(\tilde{\xi}^{2}_{i,n}\bfV^{2}_{i-1}|\bdeta)$, since the crossing terms vanish
$$
\bbE(\tilde{\xi}_{i, n}\bfV_{i-1}\tilde{\xi}_{j,n}\bfV_{j-1} | \bdeta) = \bfV_{i-1}\bfV_{j-1}\bbE(\tilde{\xi}_{i,n}\tilde{\xi}_{j, n}| \bdeta) =  \bfV_{i-1}\bfV_{j-1}\bbE(\tilde{\xi}_{i, n}| \bdeta)\bbE(\tilde{\xi}_{j, n}| \bdeta) = 0.
$$

For the first term, we have
\beqnas
& &\bbP(\|\sum^{[nt]}_{i=1}\tilde{\xi}_{i,n}\bfV_{i-1}\| \geq (\beta - K \sqrt{\delta})\sqrt{n\log n}) \\
&\leq& \frac{1}{ (\beta - K \sqrt{\delta})^{2}n\log n}\bbE\|\sum^{[nt]}_{i=1}\tilde{\xi}_{i, n}\bfV_{i-1}\|^{2}= \frac{1}{ (\beta - K \sqrt{\delta})^{2}},
\eeqnas
where the estimate of $\bbE\|\sum^{[nt]}_{i=1}\tilde{\xi}_{i, n}\bfV_{i-1}\|^{2}$ is derived from 
$$
\bbE\|\sum^{[nt]}_{i=1}\tilde{\xi}_{i, n}\bfV_{i-1}\|^{2} \leq \bbE\|\sum^{[nt]}_{i=1}\xi'_{i, n}\bfV_{i-1}\|^{2}  = O([nt]\log [nt]).
$$
Similar estimates lead to
\beqnas
\bbP(\cA^{2}_{[nt]} > K^{2}\delta n \log n) \leq \frac{1}{ K^{2}\delta n \log n}\bbE(\cA^{2}_{n}) &\leq&  \frac{1}{ K^{2}\delta}.
\eeqnas

Therefore, by taking $K \sqrt{\delta} = \frac{\beta}{2}$ we obtain
\beqnas
\bbP( \frac{1}{\sqrt{n\log n}}\sup_{t \in [0, 1]}\|\sum^{[nt]}_{i=1}\tilde{\xi}_{i, n}\bfV_{i-1}\| \geq \beta) \leq  \frac{1}{ (\beta - K \sqrt{\delta})^{2}} +  \frac{1}{ K^{2}\delta} = O(\frac{1}{\beta^{2}}), 
\eeqnas
leading to
\beqnas
\frac{1}{\sqrt{n\log n}}\bbE( \sup_{t \in [0, 1]}\|\sum^{[nt]}_{i=1}\tilde{\xi}_{i, n}\bfV_{i-1}\| ) = O(1),
\eeqnas
so that we finish our proof. 
\epf

Now we are in the position for the estimate of $I_{2}$ and $I_{3}$.

Direct computation yields
\beqnas
|\frac{1}{\sqrt{t\log t}}- \frac{1}{\bar{\xi}^{\frac{1}{2}}\sqrt{n_{t}\log n_{t}}}| = O(\frac{1}{\sqrt{t\log t}} \cdot \frac{1}{\log t}), 
\eeqnas
such that by the estimate~\eqref{est.qn} of $\bbE\|\bfQ_{n_{t}}\|$ we have 
\beqna\label{est.cont.i3}
I_{3} \leq \frac{C}{z_{1}\sqrt{t\log t} \cdot \log t}\bbE\|\bfQ_{n_{t}}\| = O(\frac{1}{z_{1}\log t}).
\eeqna

%

As a final step before the estimate of $I_{2}$, we show the tail estimate of $|\nu_{t} - n_{t}|$. To see this,  we make use of the connection between $\tau_{n}$ and $\nu_{t}$, inspired from \cite{M-T16}.  Let $\tau_{n}$ be the time up to the $n$th collision, given by
\beqnas
\tau_{n} = \sum^{n}_{i=1}\xi_{i}, \ \ \ \tau_{0} = 0
\eeqnas

First we extend $\tau_{n}$ to the continuous time. Define $\tau(u): \bbR^{+} \rightarrow \bbR^{+}$, such that
\beqnas
\tau(u) = \tau_{[u]} + \xi_{[u]+1}\{u\}. 
\eeqnas
And we  also extend $\nu_{t}$ to $\bbR^{+}$, noted as $\nu(t): \bbR^{+} \rightarrow \bbR^{+}$,
\beqnas
\nu(t) = \nu_{[t]} + \frac{\{t\}}{ \xi_{[\nu_{t}]+1}}.
\eeqnas
Then $\tau(u)$ and $\nu(t)$ are inverse mappings of each other. 

Following the proof in \cite{M-T16}, one has
\beqnas
\frac{\tau(u) - u\bar{\xi}}{\sqrt{u \log u}} \Rightarrow \cN(0, 1),
\eeqnas
as the consequence of 
\beqna\label{est.moment2.tau}
\bbE(\frac{|\tau(u) - u\bar{\xi}|^{2}}{u \log u}) = O(1).
\eeqna
Notice that $\tau \circ \nu (t) = t$, thus by taking $u = \nu(t)$, we also have
\beqna\label{est.moment2.nu}
\bbE(\frac{|t - \nu(t)\bar{\xi}|^{2}}{\nu(t) \log \nu(t)}) = O(1).
\eeqna

\blem
Let $A_{t}$ be a  positive function of $t$, then we have
\beqna\label{est.nu.n}
\bbP( |\nu_{t} - n_{t}| \geq A_{t}) \leq O(\frac{t\log t}{A^{2}_{t}}).
\eeqna
\elem
\bpf
Notice that
\beqnas
\bbP( |\nu_{t} - n_{t}| \geq A_{t}) \leq \bbP(n_{t} - \nu_{t} \geq A_{t}) + \bbP(\nu_{t} - n_{t} \geq A_{t}).
\eeqnas
For the first term, we require that
$$
c_{t} := n_{t} - A_{t} > 0,
$$
otherwise $\bbP(n_{t} - \nu_{t} \geq A_{t}) = 0$. Then we have
\beqnas
\bbP(n_{t} - \nu_{t} \geq A_{t}) 
&\leq& \bbP(\frac{n_{t} - \nu_{t}}{\sqrt{\nu_{t}\log \nu_{t}}} \geq  \frac{A_{t}}{\sqrt{c_{t}\log c_{t}}})\\
& \leq&  \frac{c_{t}\log c_{t}}{ A^{2}_{t}}\bbE\frac{|n_{t} - \nu_{t}|^{2}}{\nu_{t}\log \nu_{t}} = O(\frac{t\log t}{A^{2}_{t}}),
\eeqnas
where the last line is due to \eqref{est.moment2.nu}.

Let $m_{t} := n_{t} + A_{t}$, and we derive that
\beqnas
& &\bbP(\nu_{t} - n_{t} \geq A_{t}) =  \bbP(\tau_{m_{t}}  \leq  t) \\
&=& \bbP(\tau_{m_{t}} -   m_{t}\bar{\xi} \leq - A_{t}\bar{\xi}) \leq \bbP(|\tau_{m_{t}} -   m_{t}\bar{\xi}| \geq A_{t}\bar{\xi}) \\
&\leq& \frac{\bbE|\tau_{m_{t}} -   m_{t}\bar{\xi}|^{2}}{A^{2}_{t}\bar{\xi}^{2}} =  O(\frac{t\log t}{A^{2}_{t}}),
\eeqnas
which finishes our proof.
\epf

Now we are ready to estimate $I_{2}$. By~\eqref{est.nu.n}, we have
\beqnas
I_{2}
&\leq& \bbP( |\nu_{t} - n_{t}| \geq A_{t} ) + \bbP( \frac{1}{\sqrt{t\log t}}(\|\bfQ_{\nu_{t}} - \bfQ_{n_{t}}\|  \geq \frac{z_{1}}{3},  |\nu_{t} - n_{t}| \leq A_{t})\\
&\leq& O(\frac{t\log t}{A^{2}_{t}}) + \frac{C}{z_{1}} \frac{1}{\sqrt{t\log t}}\bbE(\|\bfQ_{\nu_{t}} - \bfQ_{n_{t}}\|, |\nu_{t} - n_{t}| \leq A_{t}).
\eeqnas
By the stationarity of $\bfQ_{n}$ and \eqref{est.qn.cont}, we have
\beqnas
& &\bbE(\|\bfQ_{\nu_{t}} - \bfQ_{n_{t}}\|, |\nu_{t} - n_{t}| \leq A_{t}) \\
&\leq& \bbE(\max_{1 \leq k \leq [A_{t}]}\|\bfQ_{n_{t} + k} - \bfQ_{n_{t}}\|, |\nu_{t} - n_{t}| \leq A_{t}) \\
&\leq& \bbE(\max_{1 \leq k \leq [A_{t}]}\|\bfQ_{k}\|) = \sqrt{A_{t}\log A_{t}},
\eeqnas
such that
\beqnas
I_{2} \leq O(\frac{t\log t}{A^{2}_{t}}) + \frac{C}{z_{1}} \frac{\sqrt{A_{t}\log A_{t}}}{\sqrt{t\log t}}.
\eeqnas
Then by taking $A_{t} =  t^{\frac{1}{2} + \epsilon}$, we have
\beqna\label{est.cont.i2}
I_{2} \leq O(\frac{|\log t|}{ t^{\epsilon}} \vee \frac{1}{z_{1}t^{\frac{1}{4} - \frac{\epsilon}{2} }}).
\eeqna

In the following we prove Theorem~\ref{thm3}. 

\bpf
Putting \eqref{est.cont.i1}, \eqref{est.cont.i2} and \eqref{est.cont.i3} into \eqref{ineq3}, we get that
\beqnas
\bbP(\|\bfW_{t} - \bfW_{n_{t}}\| \geq z_{1}) \leq  O(\frac{1}{z_{1}|\log t|}),
\eeqnas
such that by \eqref{ineq1} and \eqref{ineq2.1} (or  \eqref{ineq2.2} ), and choosing $z_{1} = \frac{1}{\sqrt{|\log t|}}$, we obtain for $d \geq 3$
\beqnas
|\bbP(\bfW_{t} \geq \bfz)  - |\bbP(\bfZ \geq \bfz)  | \leq  O(\frac{1}{(|\log t|)^{\frac{1}{4}}}),
\eeqnas
and  for $d=2$
\beqnas
|\bbP(\bfW_{t} \geq \bfz)  - |\bbP(\bfZ \geq \bfz)  | \leq  O((|\frac{\log \log t}{\log t}|)^{\frac{1}{4}}).
\eeqnas
\epf

\brmq
We point out that we are not able to prove the convergence rate for Lipschitz functions, as the analogue of Therorem~\ref{thm1} for the continuous time displacement. The main obstacle is that the term $\frac{1}{\sqrt{t \log t}}\bbE\|\bfQ_{\nu_{t}} - \bfQ_{n_{t}}\|$ is hard to control when $\nu_{t}$ is far away from $n_{t}$. 
\ermq

\section{Acknowledgements}
 I am grateful to Prof. Yong Liu, Prof. Jie Xiong for the supports. I would also like to thank Prof. Xiang-dong Li for discussions and advices. 

I would also like to express my sincere gratitude to the referee for the careful reading and helpful suggestions.



\end{document}